\documentclass[lettersize, 11pt]{article}
\usepackage{lipsum}
\usepackage{amsfonts}
\usepackage{graphicx}
\usepackage{epstopdf}
\usepackage{algorithmic}
\usepackage{color}
\usepackage{amsmath}
\usepackage{textcomp}

\usepackage{multirow}
\usepackage{algorithm}
\usepackage{algorithmic}

\usepackage{url}
\usepackage{bbm}
\usepackage{bm}
\usepackage{amsfonts}
\usepackage{verbatim} 
\usepackage{subcaption}
\usepackage{mathrsfs}
\usepackage{enumerate}

\newtheorem{theorem}{Theorem}[section]
\newtheorem{lemma}{Lemma}[section]
\newtheorem{definition}{Definition}[section]
\newtheorem{proposition}{Proposition}[section]

\newcommand{\beginproof}{\noindent{\em Proof: }}
\newcommand{\qedbox}{\hfill$\Box$ \newline\noindent}

\newtheorem{example}{Example}[section]
\newtheorem{remark}{Remark}[section]

\allowdisplaybreaks

\begin{document}

\title{ {\bf Ambiguous Chance-Constrained Binary Programs under Mean-Covariance Information}\footnote{This paper has been accepted for publication at the SIAM Journal on Optimization. Copyright \textcopyright \ 2018 Society for Industrial and Applied Mathematics}
\author{\normalsize \bf{Yiling Zhang}, \bf{Ruiwei Jiang}, and \bf{Siqian Shen} \\
{\small Department of Industrial and Operations Engineering}\\
{\small University of Michigan, Ann Arbor, MI 48109}\\[1mm]
{\small Email: \{zyiling, ruiwei, siqian\}@umich.edu}\\
}
}
\date{ }
\maketitle

\begin{abstract}
 {\color{black} We consider chance-constrained binary programs, where each row of the inequalities that involve uncertainty needs to be satisfied probabilistically.} Only the information of the mean and covariance matrix is available,  and we solve distributionally robust chance-constrained binary programs (DCBP). Using two different ambiguity sets, we equivalently reformulate the DCBPs as 0-1 second-order cone (SOC) programs. We further exploit the submodularity of 0-1 SOC constraints under special and general covariance matrices, and utilize the submodularity as well as lifting to derive extended polymatroid inequalities to strengthen the 0-1 SOC formulations. We incorporate the valid inequalities in a branch-and-cut algorithm for efficiently solving DCBPs. {\color{black} We demonstrate the computational efficacy and solution performance using diverse instances of a chance-constrained bin packing problem.}
\end{abstract}

\noindent{\it Key words:} Chance-constrained binary program; distributionally robust optimization; conic integer program; submodularity; extended polymatroid; bin packing



\section{Introduction}
\label{sec:intro}
{\color{black}
We consider chance-constrained binary programs that involve a set of individual chance constraints. More specifically, we consider $I$ individual chance constraints and, for each $i \in [I] := \{1,\ldots,I\}$, let $y_i \in \{0, 1\}^J$ be a binary decision vector such that $y_i := [y_{i1}, \ldots, y_{iJ}]^{\top}$ and $\tilde{t}_i$ be the corresponding random coefficients such that $\tilde{t}_i := [\tilde{t}_{i1}, \ldots,\tilde{t}_{iJ}]^{\top}$. Then, we consider the following individual chance constraints:
\begin{equation}
\mathbb{P} \Biggl\{ \tilde{t}_{i}^{\top} y_{i} \leq T_i\Biggr\} \ \geq \ 1 - \alpha_i \quad \forall i  \in [I], \label{cc}
\end{equation}
where $T_i \in \mathbb R$, $\mathbb{P}$ represents the joint probability distribution of $\{\tilde{t}_{ij}: i \in [I], \ j \in [J]\}$, and each $\alpha_i$ represents an allowed risk tolerance of constraint violation that often takes a small value (e.g., $\alpha_i = 0.05$).


The individual chance constraints have wide applications in service and operations management, providing an effective and convenient way of controlling capacity violation and ensuring high quality of service. For example, in surgery allocation, $y_i$ represents yes-no decisions of allocating $J$ surgeries in operating room (OR) $i$, for all $i \in [I]$. The operational time limit of each OR (i.e., $T_i$) is usually deterministic, but the processing time of each surgery (i.e., $\tilde{t}_{ij}$) is usually random due to the variety of patients, surgical teams, and surgery characteristics. Then, chance constraints \eqref{cc} make sure that each OR will not go overtime with a large probability, offering an appropriate ``end-of-the-day'' guarantee.}

Chance-constrained programs are difficult to solve, mainly because the feasible region described by constraints \eqref{cc} is non-convex in general~\cite{Prekopa2003}. Nonetheless, promising special cases have been identified to recapture the convexity of chance-constrained models. In particular, if $\{\tilde{t}_{ij}: j \in [J]\}$ are assumed to follow a Gaussian distribution with a known mean $\mu_i$ and covariance matrix $\Sigma_i$, then the chance constraints \eqref{cc} are equivalent to the second-order cone (SOC) constraints
\begin{equation}
\mu_i^{\top} y_i + \Phi^{-1}(1 - \alpha_i) \sqrt{y_i^{\top} \Sigma_i y_i} \ \leq \ T_i \quad \forall i \in [I], \label{cc-normal}
\end{equation}
where $\Phi(\cdot)$ represents the cumulative distribution function of the standard Gaussian distribution. {\color{black} In this case, feasible binary solutions $y_i$ to constraints \eqref{cc-normal} can be quickly found by off-the-shelf optimization solvers.}
In another promising research stream, the probability distribution $\mathbb{P}$ of $\tilde{t}_{ij}$ is replaced by a finite-sample approximation, leading to a sample average approximation (SAA) of the chance-constrained model~\cite{CCP_luedtke-ahmed-2008,CCPNEW_pagnoncelli2009computational}. The SAA model is then recast as a mixed-integer linear program (MILP), on which many strong valid inequalities can be derived to accelerate the branch-and-cut algorithm (see, e.g., \cite{Luedtke2010,CCP_kucukyavuz2010mixing,luedtke2014branch,song2014chance,liu2016decomposition}).

However, a basic challenge of the chance-constrained approach is that the perfect knowledge of probability distribution $\mathbb{P}$ may not be accessible. Under many circumstances, we only have a series of historical data that can be considered as samples taken from the true (while ambiguous) distribution. As a consequence, the solution obtained from a chance-constrained model can be sensitive to the choice of distribution $\mathbb{P}$ we employ in \eqref{cc} and hence perform poorly in out-of-sample tests. This phenomenon is often observed when solving stochastic programs and is called the \emph{optimizer's curse}~\cite{smith2006optimizer}. A natural way of addressing this curse is that, instead of a single estimate of $\mathbb{P}$, we employ a set of plausible probability distributions, termed the ambiguity set and denoted $\mathcal{D}$. Then, from a robust perspective, we ensure that chance constraints \eqref{cc} hold valid with regard to all probability distributions belonging to $\mathcal{D}$, i.e.,
\begin{equation}
\inf_{\mathbb{P} \in \mathcal{D}} \ \mathbb{P} \Biggl\{ \tilde{t}_{i}^{\top} y_{i} \leq T_i\Biggr\} \ \geq \ 1 - \alpha_i \quad \forall i \in [I], \label{drcc}
\end{equation}
and accordingly we call \eqref{drcc} \emph{distributionally robust chance constraints} (DRCCs).

{\color{black} In this paper, we consider distributionally robust chance-constrained binary programs (DCBP) that involve binary $y_i$-variables and DRCCs \eqref{drcc}.} Without making the Gaussian assumption on $\mathbb{P}$, we show that a DCBP is equivalent to a 0-1 SOC program when $\mathcal{D}$ is characterized by the first two moments of $\tilde{t}_{i}$. Furthermore, {\color{black} building upon existing work on valid inequalities for submodular/supermodular functions}, we exploit the submodularity of the 0-1 SOC program to generate valid inequalities. As demonstrated in extensive computational experiments of bin packing instances, these valid inequalities significantly accelerate the branch-and-cut algorithm for solving the related DCBPs. Notably, the proposed submodular approximations and the resulting valid inequalities apply to general 0-1 SOC programs than DCBPs.

The remainder of the paper is organized as follows. Section~\ref{sec:prior} reviews 
the prior work related to optimization techniques used in this paper and stochastic bin packing problems. Section~\ref{sec:model} presents two 0-1 SOC representations, respectively, for DRCCs under two moment-based ambiguity sets. Section~\ref{sec:algorithm} utilizes submodularity and lifting to generate valid inequalities to strengthen the 0-1 SOC formulations.
{\color{black} Section~\ref{sec:result} demonstrates the computational efficacy of our approaches for solving 
different 0-1 SOC reformulations of a DCBP for chance-constrained bin packing, with diverse problem sizes and parameter settings. Section~\ref{conclusion} summarizes the paper and discusses future research directions.
}

{\color{black}
\section{Prior Work}
\label{sec:prior}

Chance-constrained binary programs with uncertain technology matrix and/or right-hand side are computationally challenging, largely because (i) the non-convexity of chance constraints and (ii) the discrete variables. 
The majority of existing literature requires full distributional knowledge of the random coefficients and applies the SAA approach to approximate the models as MILPs. For example, \cite{song2014chance} considered a generic chance-constrained binary packing problem using finite samples of the random item weights, and derived lifted cover inequalities to accelerate the computation. 
For generic chance-constrained programs, the SAA approach and valid inequalities for the related MILPs have been well studied in the literature (see, e.g., \cite{CCP_luedtke-ahmed-2008,
Luedtke2010,CCP_kucukyavuz2010mixing,luedtke2014branch}). We also refer to \cite{ORA_Shylo_CC, yan2013ORP, shen-wang-CC-SS2013} for wide applications of chance-constrained binary programs, mainly in service systems and operations.
As compared to the existing work, this paper waives the assumption of full distributional information and only relies on the first two moments of the uncertainty.

Distributionally robust optimization has received growing attention, mainly because it provides effective modeling and computational approaches for handling ambiguous distributions of random variables in stochastic programming by using available distributional information.
Moment information 
has been widely used for building ambiguity sets in various distributionally robust optimization models (see, e.g., \cite{bertsimas2010models,delage2010distributionally,wiesemann2014distributionally}). Using moment-based ambiguity sets,~\cite{el2003worst, CCP_calafiore-el-ghaoui-2006, wagner2008stochastic, CCP_chen2010cvar, zymler2013distributionally, cheng2014distributionally, jiang2013data} derived exact reformulations and/or approximations for DRCCs, often in the form of semidefinite programs (SDPs). In special cases, e.g., when the first two moments are \emph{exactly} matched in the ambiguity set, the SDPs can further be simplified as SOC programs. 
While many existing ambiguity sets exactly match the first two moments of uncertainty (see, e.g., \cite{el2003worst,wagner2008stochastic,zymler2013distributionally}),~\cite{delage2010distributionally} proposed a data-driven approach to construct an ambiguity set that can model moment estimation errors. In this paper, we consider both types of moment-based ambiguity sets. To the best of our knowledge, for the first time, we provide an SOC representation of DRCCs using the general ambiguity set proposed by~\cite{delage2010distributionally}.

Meanwhile, distributionally robust optimization has received much less attention in discrete optimization problems, possibly due to the difficulty of solving 0-1 nonlinear programs. For example, most off-the-shelf solvers cannot directly handle 0-1 SDPs, which often arise from discrete optimization problems with DRCCs. To the best of our knowledge, our results on chance constraints are most related to \cite{cheng2014distributionally} that studied DRCCs in the binary knapsack problem and derived 0-1 SDP reformulations. As compared to \cite{cheng2014distributionally}, we investigate a different ambiguity set and derive a 0-1 \emph{SOC} representation. Additionally, we solve the 0-1 SOC reformulation to \emph{global} optimality instead of considering an SDP relaxation as in \cite{cheng2014distributionally}.

In the seminal work \cite{nemhauser1978analysis}, the authors identified submodularity in combinatorial and discrete optimization problems and proved a sufficient and necessary condition for 0-1 quadratic functions being submodular. We use this condition to exploit the submodularity of our 0-1 SOC reformulations. Indeed, submodular and supermodular knapsack sets (the discrete lower level set of a submodular function and discrete upper level set of a supermodular function, respectively) often arise when modeling utility, risk, and chance constraints on discrete variables. Extended polymatroid inequalities (see \cite{edmonds1970submodular,Nemhauser1999CO}) can be efficiently obtained through greedy-based separation procedures to optimize submodular/supermodular functions. Recently,~\cite{atamturk2009submodular, atamturk2015supermodular, bhardwaj2015binary} proposed cover and packing inequalities for efficiently solving submodular and supermodular knapsack sets with 0-1 variables.
Our results on valid inequalities are most related to \cite{atamturk2018network, bhardwaj2015binary} that identified a sufficient condition for the submodularity of 0-1 SOC constraints and strengthened their formulations by using the corresponding extended polymatroid inequalities (see Section 2.2 of~\cite{atamturk2018network}). In contrast, we derive a different way to exploit the submodularity of general 0-1 SOC constraints. In particular, we apply the sufficient and necessary condition derived by~\cite{nemhauser1978analysis} to search for ``optimal'' submodular approximations of the 0-1 SOC constraints (see Section~\ref{subsec:submodular}).}

The main contributions of the paper are three-fold. First, using the general moment-based ambiguity set proposed by~\cite{delage2010distributionally}, we equivalently reformulate DRCCs as 0-1 SOC constraints that can readily be solved by solvers. Second, we exploit the (hidden) submodularity of the 0-1 SOC constraints and employ extended polymatroid valid inequalities to accelerate solving DCBP. In particular, we provide an efficient way of finding ``optimal'' submodular approximations of the 0-1 SOC constraints in the original variable space, and furthermore show that any 0-1 SOC constraint possesses submodularity in a lifted space. The valid inequalities in original and lifted spaces can both be efficiently separated via the well-known greedy algorithm (see, e.g., \cite{atamturk2017polymatroid, atamturk2008polymatroids, edmonds1970submodular}). Third, we conduct extensive numerical studies to demonstrate the computational efficacy of our solution approaches.

\section{DCBP Models and Reformulations}
\label{sec:model}

We study DRCCs under two alternatives of ambiguity set $\mathcal{D}$ based on the first two moments of $\tilde{t}_{i}, \ i \in [I]$. The first ambiguity set, denoted $\mathcal{D}_1$, exactly matches the mean and covariance matrix of each $\tilde{t}_{i}$. In contrast, the second ambiguity set, denoted $\mathcal{D}_2$, considers the estimation errors of sample mean and sample covariance matrix (see \cite{delage2010distributionally}). In Section~\ref{subsec:as}, we introduce these two ambiguity sets and their calibration based on historical data. In Section~\ref{subsec:soc-ref}, we derive SOC representations of DRCC \eqref{drcc} under $\mathcal{D}_1$ and $\mathcal{D}_2$, respectively. While the former case (i.e., \eqref{drcc} under $\mathcal{D}_1$) has been well studied (see, e.g., \cite{el2003worst,CCP_calafiore-el-ghaoui-2006,zymler2013distributionally}), to the best of our knowledge, this is the first work to show the SOC representation of the latter case {\color{black}based on general covariance matrices} (i.e., \eqref{drcc} under $\mathcal{D}_2$).

\subsection{Ambiguity Sets}
\label{subsec:as}

Suppose that a series of independent historical data samples $\{\tilde{t}^n_i\}_{n = 1}^N$ are drawn from the true probability distribution $\mathbb{P}$ of $\tilde{t}_{ij}$. Then, the first two moments of $\tilde{t}_i$ can be estimated by the sample mean and sample covariance matrix
\begin{equation*}
\mu_i \ = \ \frac{1}{N} \sum_{n=1}^N \tilde{t}^n_i, \qquad \Sigma_i \ = \ \frac{1}{N} \sum_{n=1}^N (\tilde{t}^n_i - \mu_i) (\tilde{t}^n_i - \mu_i)^{\top}.
\end{equation*}
Throughout this paper, we assume that both $\Sigma_i$ and the true covariance matrix of $\tilde{t}_i$ are {\color{black}symmetric and} positive definite. As promised by the law of large numbers, as the data size $N$ grows, $\mu_i$ and $\Sigma_i$ converge to the true mean and true covariance matrix of $\tilde{t}_i$, respectively. Hence, when $N$ takes a large value, a natural choice of the ambiguity set consists of all probability distributions that match the sample moments $\mu_i$ and $\Sigma_i$, i.e.,
\begin{equation*}
\mathcal{D}_1  
\ = \ \left\{\mathbb{P} \in {\color{black}\mathcal{P}(\mathbb{R}^J)}: \
\begin{array}{l} \mathbb{E_P}[\tilde{t}_i] = \mu_i, \\[0.3cm]
\mathbb{E_P}[(\tilde{t}_i - \mu_i) (\tilde{t}_i - \mu_i)^{\top}] = \Sigma_i, \  \forall i \in [I]
\end{array}
\right\},
\end{equation*}
where {\color{black}$\mathcal{P}(\mathbb{R}^J)$} represents the set of all probability distributions on $\mathbb{R}^J$.

Under many circumstances, however, the historical data can be inadequate. With a small $N$, there may exist considerable estimation errors in $\mu_i$ and $\Sigma_i$, which brings a layer of ``moment ambiguity'' into $\mathcal{D}_1$ and adds to the existing distributional ambiguity of $\mathbb{P}$. To address the moment ambiguity and take into account the estimation errors,~\cite{delage2010distributionally} proposed an alternative ambiguity set
\begin{equation*}
\mathcal{D}_2 \ = \ \left\{\mathbb{P} \in {\color{black}\mathcal{P}(\mathbb{R}^J)}: \
\begin{array}{l} (\mathbb{E_P}[\tilde{t}_i] - \mu_i)^{\top} \Sigma_i^{-1} (\mathbb{E_P}[\tilde{t}_i] - \mu_i) \ \leq \ \gamma_1, \\[0.3cm]
\mathbb{E_P}\bigl[(\tilde{t}_i - \mu_i) (\tilde{t}_i - \mu_i)^{\top}\bigr] \ \preceq \ \gamma_2 \Sigma_i, \ \ \forall i \in [I]
\end{array}
\right\},
\end{equation*}
where $\gamma_1 > 0$ and $\gamma_2 > \max\{\gamma_1, 1\}$ represent two given parameters. Set $\mathcal{D}_2$ designates that (i) the true mean of $\tilde{t}_i$ is within an ellipsoid centered at $\mu_i$, and (ii) the true covariance matrix of $\tilde{t}_i$ is bounded from above by $\gamma_2 \Sigma_i - (\mathbb{E_P}[\tilde{t}_i] - \mu_i)(\mathbb{E_P}[\tilde{t}_i] - \mu_i)^{\top}$ (note that $\mathbb{E_P}[(\tilde{t}_i - \mu_i) (\tilde{t}_i - \mu_i)^{\top}] = \mathbb{E_P}[(\tilde{t}_i - \mathbb{E_P}[\tilde{t}_i]) (\tilde{t}_i - \mathbb{E_P}[\tilde{t}_i])^{\top}] + (\mathbb{E_P}[\tilde{t}_i] - \mu_i)(\mathbb{E_P}[\tilde{t}_i] - \mu_i)^{\top}$).

\cite{delage2010distributionally} offered a rigorous guideline for selecting $\gamma_1$ and $\gamma_2$ values (see Theorem 2 in~\cite{delage2010distributionally}) so that $\mathcal{D}_2$ includes the true distribution of $\tilde{t}_i$ with a high confidence level. In practice, we can select the values of $\gamma_1$ and $\gamma_2$ via cross validation. For example, we can divide the $N$ data points into two halves. We estimate sample moments ($\mu^1_i$, $\Sigma^1_i$) based on the first half of the data and ($\mu^2_i$, $\Sigma^2_i$) based on the second half. Then, we characterize $\mathcal{D}_2$ based on ($\mu^1_i$, $\Sigma^1_i$), and select $\gamma_1$ and $\gamma_2$ such that probability distributions with moments ($\mu^2_i$, $\Sigma^2_i$) belong to $\mathcal{D}_2$.

\subsection{SOC Representations of the DRCC} \label{subsec:soc-ref}

Now we derive SOC representations of DRCC \eqref{drcc} for all $i \in [I]$. For notation brevity, we define vector $y := [y_{i1}, \ldots, y_{iJ}]$ and omit the subscript $i$ throughout this section. First, we review the celebrated SOC representation of DRCC \eqref{drcc} under ambiguity set $\mathcal{D}_1$ in the following theorem.
\begin{theorem} \label{thm:D-1}
(Adapted from~\cite{el2003worst}, also see~\cite{wagner2008stochastic}) The DRCC \eqref{drcc} with $\mathcal{D} = \mathcal{D}_1$ is equivalent to the following SOC constraint:
\begin{equation}
\mu^{\top} y + \sqrt{\frac{1 - \alpha}{\alpha}} \sqrt{y^{\top} \Sigma y} \ \leq \ T. \label{ref-0}
\end{equation}
\end{theorem}

Theorem~\ref{thm:D-1} shows that we can recapture the convexity of DRCC~\eqref{drcc} by employing ambiguity set $\mathcal{D}_1$ to model the $\tilde{t}$ uncertainty. Perhaps more surprisingly, in this case, the convex feasible region characterized by DRCC \eqref{drcc} is SOC representable. It follows that the continuous relaxation of the DCBP model is an SOC program, {\color{black}which can be solved very efficiently by standard nonlinear optimization solvers}. 

Next, we show that DRCC \eqref{drcc} under the ambiguity set $\mathcal{D}_2$ is also SOC representable. This implies that the computational complexity of the DCBP remains the same even if we take the moment ambiguity into account. We present the main result of this section in the following theorem.
\begin{theorem} \label{thm:D-2}
DRCC \eqref{drcc} with $\mathcal{D} = \mathcal{D}_2$ is equivalent to
\begin{subequations} \label{ref-1+2}
\begin{equation}
\mu^{\top} y + \Biggl(\sqrt{\gamma_1} + \sqrt{\Bigl(\frac{1 - \alpha}{\alpha}\Bigr) (\gamma_2 - \gamma_1)}\Biggr) \sqrt{y^{\top} \Sigma y} \ \leq \ T \label{ref-1}
\end{equation}
if $\gamma_1/\gamma_2 \leq \alpha$, and is equivalent to
\begin{equation}
\mu^{\top} y + \sqrt{\frac{\gamma_2}{\alpha}} \sqrt{y^{\top} \Sigma y} \ \leq \ T \label{ref-2}
\end{equation}
\end{subequations}
if $\gamma_1/\gamma_2 > \alpha$.
\end{theorem}

{\color{black}\cite{rujeerapaiboon2015robust} considered an ambiguity set similar to $\mathcal{D}_2$ and derived an SOC representation of DRCCs in portfolio optimization under an assumption of weak sense white noise, i.e., the uncertainty is stationary and mutually uncorrelated over time (see Definition 4 and Theorem 5 in~\cite{rujeerapaiboon2015robust}). In contrast, the SOC representation in Theorem~\ref{thm:D-2} holds for general covariance matrices.} We prove Theorem~\ref{thm:D-2} in two steps. In the first step, we project the random vector $\tilde{t}$ and its ambiguity set $\mathcal{D}_2$ from $\mathbb{R}^J$ to the real line, i.e., $\mathbb{R}$. This simplifies DRCC \eqref{drcc} as involving a one-dimensional random variable. In the second step, we derive optimal (i.e., worst-case) mean and covariance matrix in $\mathcal{D}_2$ that attain the worst-case probability bound in \eqref{drcc}. We then {\color{black}apply Cantelli's inequality} to finish the representation. We present the first step of the proof in the following lemma.
\begin{lemma} \label{lemma:projection}
Let $\tilde{s}$ be a random vector in $\mathbb{R}^J$ and $\tilde{\xi}$ be a random variable in $\mathbb{R}$. For a given $y \in \mathbb{R}^J$, define ambiguity sets $\mathcal{D}_{\tilde{s}}$ and $\mathcal{D}_{\tilde{\xi}}$ as
\begin{subequations}
\begin{equation}
\mathcal{D}_{\tilde{s}} \ = \ \Bigl\{\mathbb{P} \in {\color{black}\mathcal{P}(\mathbb{R}^J)}: \ \mathbb{E_P}[\tilde{s}]^{\top}\Sigma^{-1}\mathbb{E_P}[\tilde{s}] \leq \gamma_1, \quad \mathbb{E_P}[\tilde{s} \tilde{s}^{\top}] \preceq \gamma_2 \Sigma\Bigr\} \label{D-s}
\end{equation}
and
\begin{equation}
\mathcal{D}_{\tilde{\xi}} \ = \ \Bigl\{\mathbb{P} \in {\color{black}\mathcal{P}(\mathbb{R})}: \ |\mathbb{E_P}[\tilde{\xi}]| \leq \sqrt{\gamma_1}\sqrt{y^{\top}\Sigma y}, \quad \mathbb{E_P}[\tilde{\xi}^2] \leq \gamma_2 (y^{\top}\Sigma y) \Bigr\}. \label{D-xi}
\end{equation}
\end{subequations}
Then, for any Borel measurable function $f: \mathbb{R}^J \rightarrow \mathbb{R}$, we have
\begin{equation*}
\inf_{\mathbb{P} \in \mathcal{D}_{\tilde{s}}} \mathbb{P}\{f(y^{\top}\tilde{s}) \leq 0\} \ = \ \inf_{\mathbb{P} \in \mathcal{D}_{\tilde{\xi}}} \mathbb{P}\{f(\tilde{\xi}) \leq 0\}.
\end{equation*}
\end{lemma}

\beginproof We first show that $\inf_{\mathbb{P} \in \mathcal{D}_{\tilde{s}}} \mathbb{P}\{f(y^{\top}\tilde{s}) \leq 0\} \geq \inf_{\mathbb{P} \in \mathcal{D}_{\tilde{\xi}}} \mathbb{P}\{f(\tilde{\xi}) \leq 0\}$. Pick a $\mathbb{P} \in \mathcal{D}_{\tilde{s}}$, and let $\tilde{s}$ denote the corresponding random vector and $\tilde{\xi} = y^{\top}\tilde{s}$. It follows that
\begin{subequations}
\begin{align}
\mathbb{E_P}[\tilde{\xi}] \ = & \ y^{\top} \mathbb{E_P}[\tilde{s}] \nonumber \\
\leq & \ \max_{s: \ s^{\top} \Sigma^{-1} s \ \leq \ \gamma_1} \ y^{\top} s \label{proj-note-1} \\
= & \ \max_{z: \ ||z||_2 \ \leq \sqrt{\gamma_1}} \ (\Sigma^{1/2}y)^{\top}z \ = \ \sqrt{\gamma_1} \sqrt{y^{\top} \Sigma y}, \nonumber
\end{align}
where inequality \eqref{proj-note-1} is because $\mathbb{E_P}[\tilde{s}]^{\top}\Sigma^{-1}\mathbb{E_P}[\tilde{s}] \leq \gamma_1$. Similarly, we have $\mathbb{E_P}[\tilde{\xi}] \geq -\sqrt{\gamma_1} \sqrt{y^{\top} \Sigma y}$. Meanwhile, note that
\begin{align}
\mathbb{E_P}[\tilde{\xi}^2] \ = & \ y^{\top} \mathbb{E_P}[\tilde{s}\tilde{s}^{\top}] y \nonumber \\
\leq & \ y^{\top} (\gamma_2 \Sigma) y \ = \ \gamma_2 (y^{\top} \Sigma y), \label{proj-note-2}
\end{align}
where inequality \eqref{proj-note-2} is because $\mathbb{E_P}[\tilde{s}\tilde{s}^{\top}] \preceq \gamma_2 \Sigma$. Hence, the probability distribution of $\tilde{\xi}$ belongs to $\mathcal{D}_{\tilde{\xi}}$. It follows that $\inf_{\mathbb{P} \in \mathcal{D}_{\tilde{s}}} \mathbb{P}\{f(y^{\top}\tilde{s}) \leq 0\} \geq \inf_{\mathbb{P} \in \mathcal{D}_{\tilde{\xi}}} \mathbb{P}\{f(\tilde{\xi}) \leq 0\}$.

Second, we show that $\inf_{\mathbb{P} \in \mathcal{D}_{\tilde{s}}} \mathbb{P}\{f(y^{\top}\tilde{s}) \leq 0\} \leq \inf_{\mathbb{P} \in \mathcal{D}_{\tilde{\xi}}} \mathbb{P}\{f(\tilde{\xi}) \leq 0\}$. Pick a $\mathbb{P} \in \mathcal{D}_{\tilde{\xi}}$, and let $\tilde{\xi}$ denote the corresponding random variable and $\tilde{s} = \bigl[\tilde{\xi}/(y^{\top}\Sigma y)\bigr] \Sigma y$. It follows that
\begin{align}
\mathbb{E_P}[\tilde{s}]^{\top}\Sigma^{-1}\mathbb{E_P}[\tilde{s}] \ = & \ \mathbb{E_P}[\tilde{\xi}]^2 \frac{y^{\top} \Sigma}{y^{\top}\Sigma y}\Sigma^{-1} \frac{\Sigma y}{y^{\top}\Sigma y} \ = \ \frac{\mathbb{E_P}[\tilde{\xi}]^2}{y^{\top} \Sigma y} \nonumber \\
\leq & \ \frac{\gamma_1 y^{\top} \Sigma y}{y^{\top} \Sigma y} \ = \ \gamma_1, \label{proj-note-3}
\end{align}
where inequality \eqref{proj-note-3} is because $|\mathbb{E_P}[\tilde{\xi}]| \leq \sqrt{\gamma_1} \sqrt{y^{\top} \Sigma y}$. Meanwhile, note that
\begin{align}
\mathbb{E_P}[\tilde{s} \tilde{s}^{\top}] \ = & \ \mathbb{E_P}\left[\tilde{\xi}^2\frac{\Sigma y}{y^{\top} \Sigma y}\frac{y^{\top} \Sigma}{y^{\top} \Sigma y}\right] \nonumber \\
= & \ \mathbb{E_P}[\tilde{\xi}^2] \frac{(\Sigma y)(\Sigma y)^{\top}}{(y^{\top} \Sigma y)^2} \nonumber \\
\preceq & \ \gamma_2 (y^{\top} \Sigma y) \frac{(\Sigma y)(\Sigma y)^{\top}}{(y^{\top} \Sigma y)^2} \label{proj-note-4} \\
\preceq & \ \gamma_2 (y^{\top} \Sigma y) \frac{(y^{\top} \Sigma y)\Sigma}{(y^{\top} \Sigma y)^2} \ = \ \gamma_2 \Sigma, \label{proj-note-5}
\end{align}
where inequality \eqref{proj-note-4} is because $\mathbb{E_P}[\tilde{\xi}^2] \leq \gamma_2 (y^{\top} \Sigma y)$ and inequality \eqref{proj-note-5} is because $(\Sigma y)(\Sigma y)^{\top} \preceq (y^{\top} \Sigma y) \Sigma$, which holds because for all $z \in \mathbb{R}^J$,
\begin{align*}
z^{\top} (\Sigma y)(\Sigma y)^{\top} z \ = & \ \left[ (\Sigma^{1/2}z)^{\top}(\Sigma^{1/2}y) \right]^2 \\
\leq & \ ||\Sigma^{1/2}z||^2 \ ||\Sigma^{1/2}y||^2 \ \ \ \ \ \ \ \ \mbox{(Cauchy-Schwarz inequality)} \\
= & \ (y^{\top} \Sigma y) (z^{\top} \Sigma z) \\
= & \ z^{\top} [(y^{\top} \Sigma y)\Sigma] z.
\end{align*}
\end{subequations}
Hence, the probability distribution of $\tilde{s}$ belongs to $\mathcal{D}_{\tilde{s}}$. It follows that $\inf_{\mathbb{P} \in \mathcal{D}_{\tilde{s}}} \mathbb{P}\{$ $f(y^{\top}\tilde{s}) \leq 0\} \leq \inf_{\mathbb{P} \in \mathcal{D}_{\tilde{\xi}}} \mathbb{P}\{f(\tilde{\xi}) \leq 0\}$ because $\tilde{\xi} = y^{\top} \tilde{s}$, and the proof is completed. \qedbox

\begin{remark}
\cite{popescu2007robust} and~\cite{yu2009general} showed a similar projection property for $\mathcal{D}_1$, i.e., when the first two moments of $\tilde{s}$ are exactly known. Lemma \ref{lemma:projection} employs a different transformation approach to show the projection property for $\mathcal{D}_2$ when these moments are ambiguous in the sense of~\cite{delage2010distributionally}.
\end{remark}

We are now ready to finish the proof of Theorem~\ref{thm:D-2}.\\
\noindent\emph{Proof of Theorem~\ref{thm:D-2}}: First, we define random vector $\tilde{s} = \tilde{t} - \mu$, random variable $\tilde{\xi} = y^{\top}\tilde{s}$, constant $b = T - \mu^{\top} y$, and set $S$ such that
\begin{equation*}
S \ = \ \{(\mu_1, \sigma_1) \in \mathbb{R}\times\mathbb{R}_+: \ |\mu_1| \leq \sqrt{\gamma_1}\sqrt{y^{\top} \Sigma y}, \ \mu_1^2 + \sigma_1^2 \leq \gamma_2 y^{\top} \Sigma y\}.
\end{equation*}
It follows that
\begin{subequations}
\begin{align}
\inf_{\mathbb{P} \in \mathcal{D}_2} \mathbb{P}\{\tilde{t}^{\top} y \leq T\} \ = & \ \inf_{\mathbb{P} \in \mathcal{D}_{\tilde{s}}} \mathbb{P}\{y^{\top} \tilde{s} \leq b\} \nonumber \\
= & \inf_{\mathbb{P} \in \mathcal{D}_{\tilde{\xi}}} \mathbb{P}\{\tilde{\xi} \leq b\} \label{proj-note-6} \\
= & \inf_{(\mu_1, \sigma_1) \in S} \ \inf_{\mathbb{P} \in \mathcal{D}_1(\mu_1, \sigma_1^2)} \mathbb{P}\{\tilde{\xi} \leq b\}, \label{proj-note-7}
\end{align}
where $\mathcal{D}_{\tilde{s}}$ and $\mathcal{D}_{\tilde{\xi}}$ are defined in \eqref{D-s} and \eqref{D-xi}, respectively, equality \eqref{proj-note-6} follows from Lemma \ref{lemma:projection}, and equality \eqref{proj-note-7} decomposes the optimization problem in \eqref{proj-note-6} into two layers: the outer layer searches for the optimal (i.e., worst-case) mean and covariance, while the inner layer computes the worst-case probability bound under the given mean and covariance. For the inner layer, based on {\color{black}Cantelli's inequality}, we have
\begin{equation*}
\inf_{\mathbb{P} \in \mathcal{D}_1(\mu_1, \sigma_1^2)} \mathbb{P} \{\tilde{\xi} \leq b\} = \left\{
\begin{array}{ll}
\frac{(b-\mu_1)^2}{\sigma_1^2 + (b-\mu_1)^2}, & \mbox{if $b \geq \mu_1$,} \\[0.15cm]
0, & \mbox{o.w.}
\end{array}\right.
\end{equation*}
As DRCC \eqref{drcc} states that $\inf_{\mathbb{P} \in \mathcal{D}_2} \mathbb{P}\{\tilde{t}^{\top} y \leq T\} \geq 1 - \alpha > 0$, we can assume $b \geq \mu_1$ for all $(\mu_1, \sigma_1) \in S$ without loss of generality. That is,
\begin{equation*}
b \ \geq \ \max_{(\mu_1, \sigma_1) \in S} \mu_1 \ = \ \sqrt{\gamma_1} \sqrt{y^{\top} \Sigma y}.
\end{equation*}
It follows that
\begin{align}
\inf_{\mathbb{P} \in \mathcal{D}_2} \mathbb{P}\{\tilde{t}^{\top} y \leq T\} \ = & \ \inf_{(\mu_1, \sigma_1) \in S} \ \frac{(b-\mu_1)^2}{\sigma_1^2 + (b-\mu_1)^2} \nonumber \\
= & \ \inf_{(\mu_1, \sigma_1) \in S} \ \frac{1}{\left(\frac{\sigma_1}{b-\mu_1}\right)^2 + 1}. \label{proj-note-8}
\end{align}
Note that the objective function value in \eqref{proj-note-8} decreases as $\sigma_1/(b-\mu_1)$ increases. Hence, \eqref{proj-note-8} shares optimal solutions with the following optimization problem:
\begin{equation}
\inf_{(\mu_1, \sigma_1) \in S} \ - \Bigl(\frac{\sigma_1}{b-\mu_1}\Bigr). \label{proj-note-9}
\end{equation}
The feasible region of problem \eqref{proj-note-9} is depicted in the shaded area of Figure~\ref{fig:01}.
\begin{figure}[h!]
	\centering%
    \includegraphics[width=.9\textwidth]{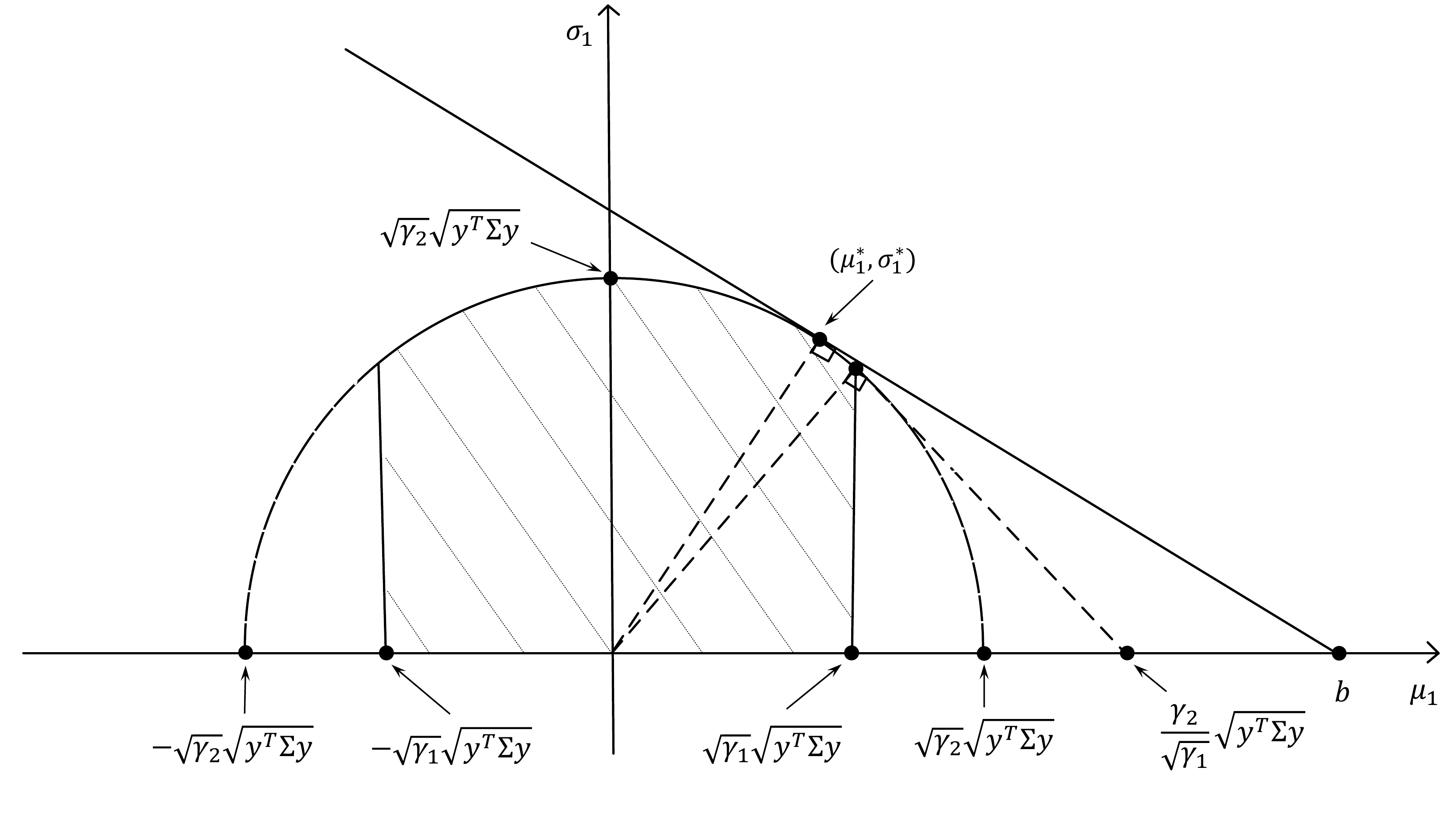}
	\caption{Graphical Solution of Problem \eqref{proj-note-9}}
	\label{fig:01}
\end{figure}
Furthermore, we note that the objective function of \eqref{proj-note-9} equals to the slope of the straight line connecting points $(b, 0)$ and $(\mu_1, \sigma_1)$ (see Figure~\ref{fig:01} for an example). It follows that an optimal solution ($\mu_1^*, \sigma_1^*$) to problem \eqref{proj-note-9}, and so to problem \eqref{proj-note-8}, lies in one of the following two cases:
\begin{description}
\item[Case 1.] If $\sqrt{\gamma_1}\sqrt{y^{\top} \Sigma y} \leq b \leq (\gamma_2/\sqrt{\gamma_1}) \sqrt{y^{\top} \Sigma y}$, then $\mu_1^* = \sqrt{\gamma_1}\sqrt{y^{\top} \Sigma y}$ and $\sigma_1^* = \sqrt{\gamma_2 - \gamma_1}\sqrt{y^{\top} \Sigma y}$.
\item[Case 2.] If $b > (\gamma_2/\sqrt{\gamma_1}) \sqrt{y^{\top} \Sigma y}$, then $\mu_1^* = (\gamma_2 y^{\top} \Sigma y)/b$ and \\$\sigma_1^* = \sqrt{\gamma_2 y^{\top} \Sigma y - (\gamma_2 y^{\top} \Sigma y)^2/b^2}$.
\end{description}
Denoting $\kappa (b, y) = \frac{b}{\sqrt{y^{\top} \Sigma y}}$, we have
\begin{equation}
\inf_{\mathbb{P} \in \mathcal{D}_2} \mathbb{P}\{\tilde{t}^{\top} y \leq T\} \ = \ \left\{
\begin{array}{ll} \frac{1}{\left( \frac{\sqrt{\gamma_2 - \gamma_1}}{\kappa(b, y) - \sqrt{\gamma_1}} \right)^2 + 1}, & \mbox{if } \sqrt{\gamma_1} \leq \kappa(b, y) \leq \frac{\gamma_2}{\sqrt{\gamma_1}}, \\[0.75cm]
{\color{black}\frac{\kappa(b,y)^2 - \gamma_2}{\kappa(b,y)^2}}, & \mbox{if } \kappa(b, y) > \frac{\gamma_2}{\sqrt{\gamma_1}}.
\end{array}
\right. \label{wc-prob-bound}
\end{equation}

Second, based on \eqref{wc-prob-bound}, the DRCC $\inf_{\mathbb{P} \in \mathcal{D}_2} \mathbb{P}\{\tilde{t}^{\top} y \leq T\} \geq 1 - \alpha$ has the following representations:
\begin{equation*}
\mbox{DRCC} \ \Leftrightarrow \ \left\{
\begin{array}{ll} \frac{b}{\sqrt{y^{\top} \Sigma y}} \ \geq \ \sqrt{\gamma_1} + \sqrt{\Bigl(\frac{1 - \alpha}{\alpha}\Bigr)(\gamma_2 - \gamma_1)}, & \mbox{if } \sqrt{\gamma_1} \leq \frac{b}{\sqrt{y^{\top} \Sigma y}} \leq \frac{\gamma_2}{\sqrt{\gamma_1}}, \\[0.2cm]
\frac{b}{\sqrt{y^{\top} \Sigma y}} \ \geq \ \sqrt{\frac{\gamma_2}{\alpha}}, & \mbox{if } \frac{b}{\sqrt{y^{\top} \Sigma y}} > \frac{\gamma_2}{\sqrt{\gamma_1}}.
\end{array}
\right.
\end{equation*}
\end{subequations}

It follows that DRCC \eqref{drcc} is equivalent to an SOC constraint by discussing the following two cases:
\begin{description}
\item[Case 1.] If $\gamma_1 / \gamma_2 \leq \alpha$, then $\gamma_2/\sqrt{\gamma_1} \geq \sqrt{\gamma_1} + \sqrt{[(1 - \alpha)/\alpha](\gamma_2 - \gamma_1)}$ and $\gamma_2/\sqrt{\gamma_1} \geq \sqrt{\gamma_2/\alpha}$. It follows that {\color{black} (i) if $\sqrt{\gamma_1} \leq b/\sqrt{y^{\top} \Sigma y} \leq \gamma_2/\sqrt{\gamma_1}$, then DRCC is equivalent to $b/\sqrt{y^{\top} \Sigma y} \geq \sqrt{\gamma_1} + \sqrt{[(1 - \alpha)/\alpha](\gamma_2 - \gamma_1)}$ and (ii) if $b/\sqrt{y^{\top} \Sigma y}$ $> \gamma_2/\sqrt{\gamma_1}$, then DRCC always holds. Combining sub-cases (i) and (ii) yields that} DRCC \eqref{drcc} is equivalent to $b/\sqrt{y^{\top} \Sigma y} \geq \sqrt{\gamma_1} + \sqrt{[(1 - \alpha)/\alpha](\gamma_2 - \gamma_1)}$.
\item[Case 2.] If $\gamma_1 / \gamma_2 > \alpha$, then $\gamma_2/\sqrt{\gamma_1} < \sqrt{\gamma_1} + \sqrt{[(1 - \alpha)/\alpha](\gamma_2 - \gamma_1)}$ and $\gamma_2/\sqrt{\gamma_1} < \sqrt{\gamma_2/\alpha}$. It follows that {\color{black} (i) if $\sqrt{\gamma_1} \leq b/\sqrt{y^{\top} \Sigma y} \leq \gamma_2/\sqrt{\gamma_1}$, then DRCC always fails and (ii) if $b/\sqrt{y^{\top} \Sigma y} > \gamma_2/\sqrt{\gamma_1}$, then DRCC is equivalent to $b/\sqrt{y^{\top} \Sigma y} \geq \sqrt{\gamma_2/\alpha}$. Combining sub-cases (i) and (ii) yields that} DRCC \eqref{drcc} is equivalent to $b/\sqrt{y^{\top} \Sigma y} \geq \sqrt{\gamma_2/\alpha}$.
\end{description}
The proofs of the above two cases are both completed by the definition of $b$.
\qedbox

	\begin{figure}[hbtp]
		\centering
      \includegraphics[width =0.6\textwidth]{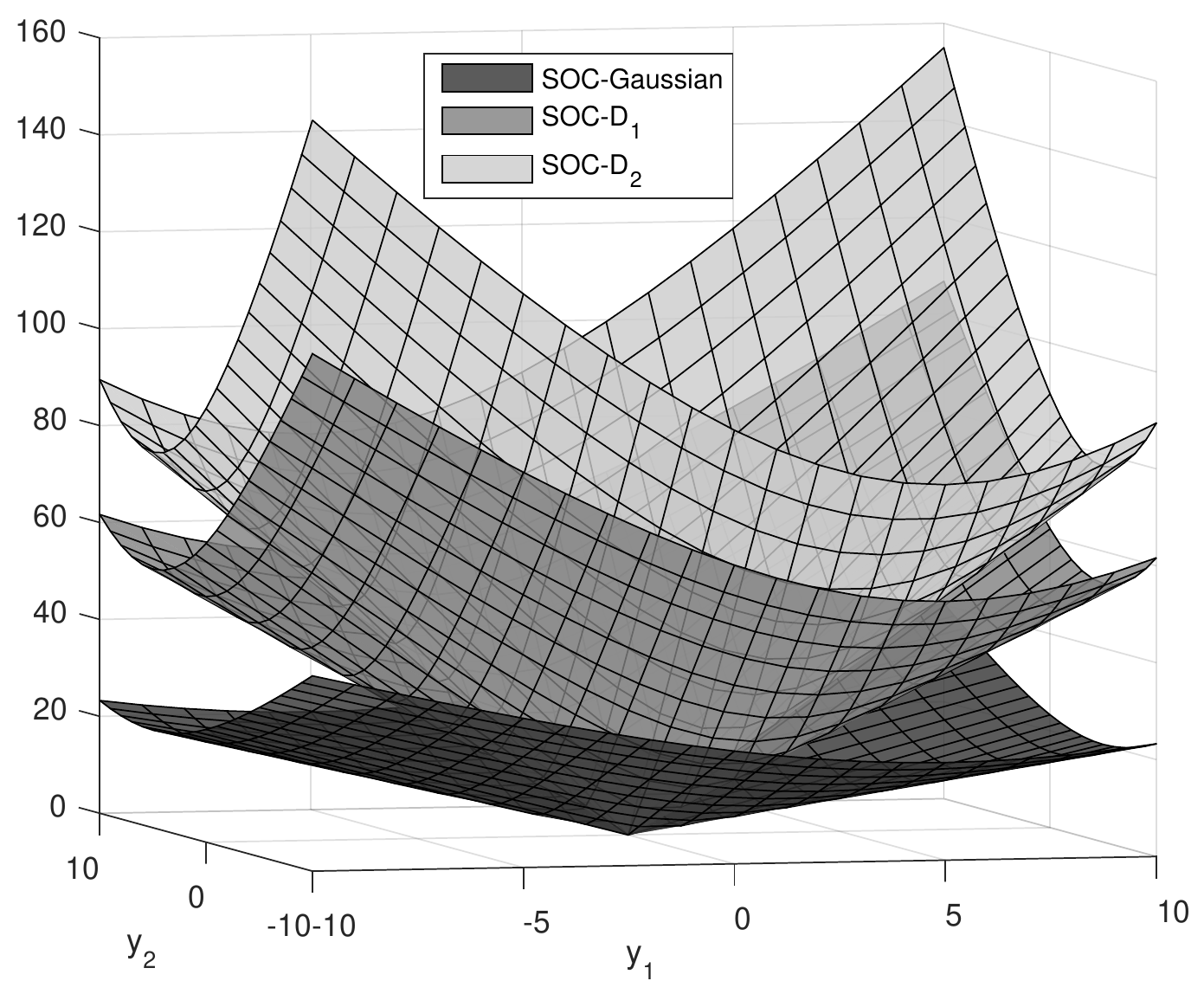}
		\caption{Illustration of the Three SOC Reformulations \eqref{cc-normal},   \eqref{ref-0}, and \eqref{ref-2} of DRCC \eqref{drcc}}
		\label{fig:threeOmega}
	\end{figure}

To sum up, we have  two exact 0-1 SOC constraint reformulations of DRCC \eqref{drcc} under ambiguity sets $\mathcal{D}_1$ and $\mathcal{D}_2$, being constraints \eqref{ref-0}, \eqref{ref-1}/\eqref{ref-2}, respectively.

{\color{black}
\begin{example}
\label{ex:threeOmega}
We consider a DRCC \eqref{drcc} with $I = 1$ (the subscript $i$ is hence omitted), $J = 2$, $1- \alpha = 95\%$, mean vector $\mu = [0 \ 0]^\top$ and covariance matrix
$
\Sigma \ = \ \begin{bmatrix}
2 & 1 \\
1 & 2
\end{bmatrix}.$
For the ambiguity set $\mathcal D_2$, we set $\gamma_1 = 1$ and $\gamma_2=2$. We note that $\Phi^{-1}(1 - \alpha)  = 1.6449$ in the SOC reformulation \eqref{cc-normal}, $\sqrt{(1-\alpha)/\alpha} = 4.3589$ in the SOC reformulation \eqref{ref-0}, and $\sqrt{\gamma_2/\alpha} = 6.3246$ in the SOC reformulation \eqref{ref-2} (since $\gamma_1/\gamma_2 > \alpha$). Without specifying the value of $T$, we depict the boundaries of the second-order cones associated with the three SOC reformulations in Figure~\ref{fig:threeOmega}. From this figure, we observe that the Gaussian approximation leads to the widest cone and so the largest SOC feasible region corresponding to DRCC \eqref{drcc}.
 \qedbox
\end{example}
}

\section{Valid Inequalities for DCBP}
\label{sec:algorithm}

Although 0-1 SOC constraint reformulations can be directly handled by the off-the-shelf solvers, as we report in Section \ref{sec:result}, {\color{black} the resultant 0-1 SOC programs} are often time-consuming to solve, primarily because of the binary restrictions on variables. In this section, we derive valid inequalities for DRCC \eqref{drcc}, with the objective of accelerating the branch-and-cut algorithm for solving DCBP with individual DRCCs and also general 0-1 SOC programs in commercial solvers. 
Specifically, we exploit the submodularity of 0-1 SOC constraints.
 In Section \ref{subsec:submodular}, we derive a sufficient condition for submodularity and two approximations of the 0-1 SOC constraints that satisfy this condition. Using the submodular approximations, we derive extended polymatroid inequalities. In Section \ref{subsec:bp}, we show the submodularity of the 0-1 SOC constraints in a lifted (i.e., higher-dimensional) space.
 While the extended polymatroid inequalities are well-known (see, e.g., \cite{atamturk2017polymatroid,atamturk2008polymatroids,edmonds1970submodular}), to the best of our knowledge, the two submodular approximations and the submodularity of 0-1 SOC constraints in the lifted space appear to be new and have not been studied in any literature.

\subsection{Submodularity of the 0-1 SOC Constraints} \label{subsec:submodular}
We consider SOC constraints of the form
\begin{equation}
\mu^{\top} y + \sqrt{y^{\top} \Lambda y} \ \leq \ T, \label{soc}
\end{equation}
where $\Lambda$ represents a $J\times J$ positive definite matrix. Note that all SOC reformulations \eqref{ref-0}, \eqref{ref-1}, and \eqref{ref-2} derived in Section \ref{sec:model}, as well as the SOC reformulation \eqref{cc-normal} of a chance constraint with Gaussian uncertainty, possess the form of \eqref{soc}. Before investigating the submodularity of \eqref{soc}, we review the definitions of submodular functions and extended polymatroids.
\begin{definition}\label{def:sub}
Define the collection of set $[J]$'s subsets $\mathcal{C} := \{R: R \subseteq [J]\}$. A set function $h$: $\mathcal{C} \rightarrow \mathbb R$ is called submodular if and only if $h(R \cup \{j\}) - h(R) \geq h(S \cup \{j\}) - h(S)$ for all $R \subseteq S \subseteq [J]$ and all $j \in [J]\setminus S$.
\end{definition}
Throughout this section, we refer to a set function as $h(R)$ and $h(y)$ interchangeably, where $y \in \{0, 1\}^J$ represents the indicating vector for subset $R \subseteq [J]$, i.e., $y_j = 1$ if $j \in R$ and $y_j = 0$ otherwise.
\begin{definition}
For a submodular function $h(y)$, the polyhedron
\begin{equation*}
\mbox{EP}_h = \{\pi \in \mathbb R^J: \ \pi(R) \leq h(R), \quad \forall R \subseteq [J]\}
\end{equation*}
is called an extended polymatroid associated with $h(y)$, where $\pi(R) = \sum_{j \in R} \pi_j$.
\end{definition}

For a submodular function $h(y)$ with $h(\emptyset) = 0$, inequality
\begin{equation}
\pi^{\top} y \leq t, \label{ep}
\end{equation}
termed an \emph{extended polymatroid inequality}, is valid for the epigraph of $h$, i.e.,
$$\{(y, t) \in \{0, 1\}^J \times \mathbb{R}: t \geq h(y)\}, \mbox{ if and only if $\pi \in \mbox{EP}_h$ {\color{black} (see \cite{edmonds1970submodular})} }.$$
Furthermore, the separation of \eqref{ep} can be efficiently done by a greedy algorithm~\cite{edmonds1970submodular}, which we briefly describe in Algorithm \ref{algo:greedy}.
\begin{algorithm}[ht!]
\caption{A greedy algorithm for separating extended polymatroid inequalities \eqref{ep}.}
\label{algo:greedy}
    \begin{algorithmic}[1]
\REQUIRE A point $(\hat{y}, \hat{t})$ with $\hat{y} \in [0, 1]^J$, $\hat{t} \in \mathbb{R}$, and a function $h$ that is submodular on $\{0,1\}^J$. 
\STATE Sort the entries in $\hat{y}$ such that $y_{(1)} \geq \cdots \geq y_{(J)}$. Obtain the permutation $\{(1), \ldots, (J)\}$ of $[J]$.
\STATE Letting $R_{(j)} := \{(1), \ldots, (j)\}, \ \forall j \in \mathcal{J}$, compute $\hat{\pi}_{(1)} = h(R_{(1)})$, and $\hat{\pi}_{(j)} = h(R_{(j)}) - h(R_{(j-1)})$ for $j = 2, \ldots, J$.
\IF{$\hat{\pi}^{\top} \hat{y} > \hat{t}$}
\STATE We generate a valid extended polymatroid inequality $\hat{\pi}^{\top} y \leq t$.
\ELSE
\STATE The current solution $(\hat{y}, \hat{t})$ satisfies $h(\hat{y}) \leq \hat{t}$.
\ENDIF
\RETURN either $(\hat{y}, \hat{t})$ is feasible, or a violated extended polymatroid inequality $\hat{\pi}^{\top} y \leq t$.
        \end{algorithmic}
\end{algorithm}

%

The strength and efficient separation of the extended polymatroid inequality motivate us to explore the submodularity of function $g(y) := \mu^{\top} y + \sqrt{y^{\top} \Lambda y}$. In a special case, $\Lambda$ is assumed to be a diagonal matrix and so the random coefficients $\tilde{t}_{ij}, j \in [J]$ for the same $i$ are uncorrelated. In this case,~\cite{atamturk2008polymatroids} successfully show that $g(y)$ is submodular. As a result, we can strengthen a DCBP by incorporating extended polymatroid inequalities in the form $\pi^{\top} y \leq T$, where $\pi \in \mbox{EP}_g$. Unfortunately, the submodularity of $g(y)$ quickly fades when the off-diagonal entries of $\Lambda$ become non-zero, e.g., when $\Lambda$ is associated with a general covariance matrix. We present an example as follows.
\begin{example}
\label{ex:general-cov}
Suppose that $[J] = \{1, 2, 3\}$, $\mu = [0, 0, 0]^\top$, and
\begin{equation*}
\Lambda \ = \ \begin{bmatrix}
0.6 & -0.2 & 0.2 \\
-0.2 & 0.7 & 0.1 \\
0.2 & 0.1 & 0.6
\end{bmatrix}.
\end{equation*}
The three eigenvalues of $\Lambda$ are 0.2881, 0.7432, and 0.8687, and so $\Lambda \succ 0$. However, function $g(y) = \mu^{\top} y + \sqrt{y^{\top} \Lambda y}$ is not submodular because $g(R \cup \{j\}) - g(R) < g(S \cup \{j\}) - g(S)$, where $R = \{1\}$, $S = \{1, 2\}$, and $j = 3$. \qedbox
\end{example}

In this section, we provide a sufficient condition for function $g(y)$ being submodular for general $\Lambda$. To this end, we apply a necessary and sufficient condition derived in~\cite{nemhauser1978analysis} for quadratic function $y^{\top}\Lambda y$ being submodular (see the second paragraph on Page 276 of~\cite{nemhauser1978analysis}, following Proposition 3.5). We summarize this condition in the following theorem.
\begin{theorem} \label{thm:submodular}
(\cite{nemhauser1978analysis}) Define function $h: \{0, 1\}^J \rightarrow \mathbb{R}$ such that $h(y) := y^{\top}\Lambda y$, where $\Lambda \in \mathbb{R}^{J\times J}$ represents a symmetric matrix. Then, $h(y)$ is submodular if and only if $\Lambda_{rs} \leq 0$ for all $r, s \in [J]$ and $r \neq s$.
\end{theorem}

Note that Theorem~\ref{thm:submodular} does not assume $\Lambda \succeq 0$ and so can be applied to general (convex or non-convex) quadratic functions. Theorem~\ref{thm:submodular} leads to a sufficient condition for function $g(y)$ being submodular, {\color{black}as presented in the following proposition.}
\begin{proposition} \label{prop-sufficient}
Let $\Lambda \in \mathbb{R}^{J \times J}$ represent a symmetric and positive semidefinite matrix that satisfies (i) $2\sum_{s=1}^J \Lambda_{rs} \geq \Lambda_{rr}$ for all $r \in [J]$ and (ii) $\Lambda_{rs} \leq 0$ for all $r, s \in [J]$ and $r \neq s$. Then, function $g(y) = \mu^{\top} y + \sqrt{y^{\top} \Lambda y}$ is submodular.
\end{proposition}
\beginproof As $\mu^{\top}y$ is submodular in $y$, it suffices to prove that $\sqrt{y^{\top} \Lambda y}$ is submodular. Hence, we can assume $\mu = 0$ without loss of generality. We let $f(x) = \sqrt{x}$ and $h(y) = y^{\top} \Lambda y$. Then, $g(y) = f(h(y))$.

First, we note that $h(y \vee e_r) \geq h(y)$ for all $y \in \{0, 1\}^J$ and $r \in [J]$, where $a \vee b = [\max\{a_1, b_1\}, \ldots, \max\{a_J, b_J\}]^{\top}$ for $a, b \in \mathbb{R}^J$. Indeed, if $y_r = 1$ then $y \vee e_r = y$ and so $h(y \vee e_r) = h(y)$. Otherwise, if $y_r = 0$, then $y \vee e_r = y + e_r$ and so
\begin{align*}
h(y \vee e_r) \ = & \ y^{\top}\Lambda y + 2e_r^{\top}\Lambda y + e_r^{\top}\Lambda e_r \\
= & \ y^{\top}\Lambda y + 2\sum_{s: \ y_s = 1} \Lambda_{rs} + \Lambda_{rr}\\
\geq & \ y^{\top}\Lambda y + 2\sum_{\substack{s = 1 \\ s \neq r}}^J \Lambda_{rs} + \Lambda_{rr}\\
\geq & \ y^{\top}\Lambda y \ = \ h(y),
\end{align*}
where the first inequality is due to $y_r = 0$ and condition (ii), and the second inequality is due to condition (i). It follows that $h(y') \geq h(y)$ for all $y, y' \in \{0, 1\}^J$ such that $y' \geq y$. Hence, $h(y)$ is increasing.

Second, based on Theorem~\ref{thm:submodular}, $h(y)$ is submodular due to condition (ii). It follows that $g(y) = f(h(y))$ is submodular because function $f$ is concave and nondecreasing and function $h$ is submodular and increasing (see, e.g., Proposition 2.2.5 in \cite{simchi2014logic})\footnote{Proposition 2.2.5 in~\cite{simchi2014logic} assumes that $y \in \mathbb{R}^n$ and provides a sufficient condition for $g$ being supermodular. It can be shown that a similar proof of this proposition applies to our case.}. \qedbox

Proposition~\ref{prop-sufficient} generalizes the sufficient condition in~\cite{atamturk2008polymatroids} because conditions (i)--(ii) are automatically satisfied if $\Lambda$ is diagonal and positive definite.

For general $\Lambda \succ 0$ that does not satisfy sufficient conditions (i)--(ii), we can approximate SOC constraint \eqref{soc} by replacing $\Lambda$ with a matrix $\Delta$ that satisfies these conditions. We derive relaxed and conservative submodular approximations of constraint \eqref{soc} in the following {\color{black}proposition}.
\begin{proposition} \label{prop-approx}
Constraint \eqref{soc} implies the SOC constraint
\begin{equation}
\mu^{\top} y + \sqrt{y^{\top} \Delta^{\mbox{\tiny L}} y} \ \leq \ T, \label{soc-relaxed}
\end{equation}
where function $g^{\mbox{\tiny L}}(y) := \mu^{\top} y + \sqrt{y^{\top} \Delta^{\mbox{\tiny L}} y}$ is submodular and $\Delta^{\mbox{\tiny L}}$ is an optimal solution of SDP
\begin{subequations}
\label{eq:relaxSDP}
\begin{align}
\min_{\Delta} \ & \ ||\Delta - \Lambda||_2 \label{relax-obj} \\
\mbox{s.t.} \ & \ 0 \ \preceq \ \Delta \ \preceq \ \Lambda, \label{relax-con-1} \\
& \ 2\sum_{s=1}^J \Delta_{rs} \geq \Delta_{rr}, \ \ \forall r \in [J], \label{relax-con-2} \\
& \ \Delta_{rs} \leq 0, \ \ \forall r, s \in [J] \mbox{ and } r \neq s. \label{relax-con-3}
\end{align}
\end{subequations}
Additionally, constraint \eqref{soc} is implied by the SOC constraint
\begin{equation}
\mu^{\top} y + \sqrt{y^{\top} \Delta^{\mbox{\tiny U}} y} \ \leq \ T, \label{soc-conservative}
\end{equation}
where function $g^{\mbox{\tiny U}}(y) := \mu^{\top} y + \sqrt{y^{\top} \Delta^{\mbox{\tiny U}} y}$ is submodular and $\Delta^{\mbox{\tiny U}}$ is an optimal solution of SDP
\begin{subequations}
\label{eq:conservativeSDP}
\begin{align}
\min_{\Delta} \ & \ ||\Delta - \Lambda||_2 \label{conservative-obj} \\
\mbox{s.t.} \ & \ \Delta \ \succeq \ \Lambda, \ \ \mbox{\eqref{relax-con-2}--\eqref{relax-con-3}}. \label{conservative-con-1}
\end{align}
\end{subequations}
\end{proposition}
\beginproof By construction, $g^{\mbox{\tiny L}}(y)$ is submodular because $\Delta^{\mbox{\tiny L}}$ satisfies constraints \eqref{relax-con-2}--\eqref{relax-con-3} and so conditions (i)--(ii). Additionally, constraint \eqref{soc} implies \eqref{soc-relaxed} because $\Delta^{\mbox{\tiny L}}$ satisfies constraint \eqref{relax-con-1} and so $\Delta^{\mbox{\tiny L}} \preceq \Lambda$. Similarly, we obtain that $g^{\mbox{\tiny U}}(y)$ is submodular and constraint \eqref{soc} is implied by \eqref{soc-conservative}. \qedbox

Note that there always exist matrices $\Delta^{\mbox{\tiny L}}$ and $\Delta^{\mbox{\tiny U}}$ that are feasible to SDPs \eqref{relax-obj}--\eqref{relax-con-3} and \eqref{conservative-obj}--\eqref{conservative-con-1}, respectively. For example, $\mbox{diag}(\lambda_{\mbox{\tiny min}}, \ldots, \lambda_{\mbox{\tiny min}}) \in \mathbb{R}^{J\times J}$ satisfy constraints \eqref{relax-con-1}--\eqref{relax-con-3}, where $\lambda_{\mbox{\tiny min}}$ represents the smallest eigenvalue of matrix $\Lambda$. Additionally, $\mbox{diag}(\lambda_{\mbox{\tiny max}}, \ldots, \lambda_{\mbox{\tiny max}}) \in \mathbb{R}^{J\times J}$ satisfy constraints \eqref{conservative-con-1}, where $\lambda_{\mbox{\tiny max}}$ represents the largest eigenvalue of matrix $\Lambda$.

By minimizing the $\ell^2$ distance between $\Delta$ and $\Lambda$ in objective functions \eqref{relax-obj} and \eqref{conservative-obj}, we find ``optimal'' approximations of $\Lambda$ that satisfies sufficient conditions (i)--(ii) in Proposition~\ref{prop-sufficient}. Accordingly, we obtain ``optimal'' submodular approximations of the 0-1 SOC constraint \eqref{soc}. There are many possible alternatives of the $\ell^2$ norm in \eqref{relax-obj} and \eqref{conservative-obj}. For example, formulations \eqref{relax-obj}--\eqref{relax-con-3} and \eqref{conservative-obj}--\eqref{conservative-con-1} remain SDPs if the $\ell^2$ norm is replaced by the $\ell^1$ norm or the $\ell^{\infty}$ norm. {\color{black} We have empirically tested $\ell^1$, $\ell^2$, and $\ell^{\infty}$ norms based on a server allocation problem (see Section \ref{sec:comp_setup} for a brief description) and the $\ell^2$ norm leads to the largest improvement on CPU time.} In computation, we only need to solve these two SDPs once to obtain $\Delta^{\mbox{\tiny L}}$ and $\Delta^{\mbox{\tiny U}}$. Then, extended polymatroid inequalities can be obtained from the relaxed approximation \eqref{soc-relaxed}. Additionally, the conservative approximation \eqref{soc-conservative} leads to an upper bound of the optimal objective value of the related DCBP.

{\color{black}
\begin{example}
\label{ex:ApproxMatrices}
Recall the $3\times 3$ matrix $\Lambda$ in Example~\ref{ex:general-cov} and the corresponding function $g(y) = \mu^{\top} y +  \sqrt{0.6 y^2_1 + 0.7 y^2_2 + 0.6 y^2_3 - 0.4 y_1 y_2 + 0.4 y_1 y_3 +  0.2 y_2 y_3}$ being not submodular. We set $\mu = [0, 0, 0]^\top$ and apply Proposition~\ref{prop-approx} to optimize the two SDPs \eqref{eq:relaxSDP} and \eqref{eq:conservativeSDP}, yielding the following two positive semidefinite matrices:
\begin{equation*}
\Delta^{\mbox{\tiny L}}  \ = \ \begin{bmatrix}
0.35 & -0.15 & 0\\
 -0.15 & 0.37 & 0\\
 0 & 0 & 0.38
\end{bmatrix} \mbox{ and}
\end{equation*}
\begin{equation*}
\Delta^{\mbox{\tiny U}} \ = \  \begin{bmatrix}
0.83 & -0.22 & 0\\
 -0.22 & 0.95 & 0\\
  0 & 0 & 0.82
\end{bmatrix}.
\end{equation*}
Due to Proposition~\ref{prop-sufficient}, $g^{\mbox{\tiny L}}(y) := \mu^{\top} y + \sqrt{y^{\top} \Delta^{\mbox{\tiny L}} y } = \sqrt{0.35 y^2_1  + 0.37 y^2_2 + 0.38 y^2_3 - 0.3 y_1 y_2} $ and $g^{\mbox{\tiny U}}(y) := \mu^{\top} y + \sqrt{y^{\top} \Delta^{\mbox{\tiny U}} y }= $
$\sqrt{0.83y^2_1 + 0.95 y^2_2 + 0.82 y^2_3 - 0.44 y_1 y_2}$ are submodular.

Now suppose that $T = 0.8$ and we are given $\hat{y} = [\hat{y}_1, \hat{y}_2, \hat{y}_3]^{\top} =  [1,0.5,0.9]^{\top}$ with $g(\hat{y}) =  \mu^{\top} \hat{y} + \sqrt{\hat{y}^{\top}\Lambda \hat{y}} = 1.23 > 0.8 = T$. First, with respect to constraint $g^{\mbox{\tiny L}}(y) \leq T$, we follow Algorithm~\ref{algo:greedy} to find an extended polymatroid inequality and note that this inequality is also valid for constraint $g(y) \leq T$. Specifically, we sort the entries of $\hat{y}$ to obtain $\hat{y}_1 \geq \hat{y}_3 \geq \hat{y}_2$ and $\{1, 3, 2\}$, the corresponding permutation. It follows that $R_{(1)} = \{1\}$, $R_{(2)} = \{1,3\}$, and $R_{(3)} = \{1, 3, 2\}$. Hence, $\hat{\pi}_{(1)} = g^{\mbox{\tiny L}}([1, 0, 0]^{\top}) = 0.59$, $\hat{\pi}_{(2)} = g^{\mbox{\tiny L}}([1, 0, 1]^{\top}) - g^{\mbox{\tiny L}}([1, 0, 0]^{\top}) = 0.26$, and $\hat{\pi}_{(3)}  = g^{\mbox{\tiny L}}([1, 1, 1]^{\top}) - g^{\mbox{\tiny L}}([1, 0, 1]^{\top}) = 0.04$. This generates the extended polymatroid inequality $0.59 y_1 + 0.26 y_3 + 0.04 y_2 \leq 0.8$ that cuts off $\hat{y}$. Second, we can replace constraint $g(y)\leq T$ with $g^{\mbox{\tiny U}}(y) \leq T$ in a DCBP to obtain a conservative approximation.
 \qedbox
\end{example}
}

\subsection{Valid Inequalities in a Lifted Space} \label{subsec:bp}

In Section \ref{subsec:submodular}, we derive extended polymatroid inequalities based on a relaxed approximation of SOC constraint \eqref{soc}.
In this section, we show that the submodularity of \eqref{soc} holds for general $\Lambda$ in a lifted (i.e., higher-dimensional) space. Accordingly, we derive extended polymatroid inequalities in the lifted space. {\color{black}We note that this approach was first proposed in~\cite{bhardwaj2015binary} to test quadratic constrained problems (see Section 3.6.2 in~\cite{bhardwaj2015binary}).}

To this end, we reformulate constraint \eqref{soc} as $\mu^{\top} y \leq T$ and $y^{\top} \Lambda y \leq (T - \mu^{\top} y)^2$, i.e., $y^{\top} (\Lambda - \mu\mu^{\top}) y + 2T\mu^{\top}y \leq T^2$. {\color{black}Note that $\mu^{\top} y \leq T$ is because $T - \mu^{\top} y \geq \sqrt{y^{\top}\Lambda y} \geq 0$ by \eqref{soc}.} Then, we define $w_{jk} = y_j y_k$ for all $j, k  \in [J]$ and augment vector $y$ to vector $v = [y_1, \ldots, y_J, w_{11}, \ldots, w_{1J}, w_{21}, \ldots, w_{JJ}]^{\top}$. We can incorporate the following McCormick inequalities to define each $w_{ij}$:
\begin{equation}
w_{jk} \leq y_j, \quad w_{jk} \leq y_k, \quad w_{jk} \geq y_j + y_k - 1, \quad w_{ij} \geq 0. \label{mccormick}
\end{equation}
Accordingly, we rewrite \eqref{soc} as $[\mu^{\top}, 0^{\top}] v \leq T$ and
\begin{equation}
\label{eq:lifted:v}
a^{\top} v + v^{\top} B_{\mbox{\tiny N}} v \leq T^2,
\end{equation}
where we decompose $(\Lambda - \mu\mu^{\top})$ to be the sum of two matrices, one containing all positive entries and the other containing all nonpositive entries. Accordingly, we define vector $a := [2T\mu; B_{\mbox{\tiny P}}]^{\top}$ with $B_{\mbox{\tiny P}} \in \mathbb{R}^{J^2}_+$ representing all the positive entries after vectorization, and matrix $B_{\mbox{\tiny N}} \in \mathbb{R}^{(J + J^2)\times (J + J^2)}_-$ collects all nonpositive entries, i.e., {\color{black}$B_{\mbox{\tiny N}} = \begin{bmatrix} -(\Lambda - \mu\mu^{\top})_- & 0 \\ 0 & 0 \end{bmatrix}$}, where $(x)_- = - \min\{0, x\}$ for $x \in \mathbb{R}$.

As $a^{\top} v + v^{\top} B_{\mbox{\tiny N}} v$ is a submodular function of $v$ by Theorem~\ref{thm:submodular}, we can incorporate extended polymatroid inequalities to strengthen the lifted SOC constraints \eqref{soc}. We summarize this result in the following {\color{black}proposition}.
\begin{proposition}
\label{theorem:lifted}
{\color{black}(See also Section 3.6.2 in~\cite{bhardwaj2015binary})} Define function $h:\{0, 1\}^{J + J^2} \rightarrow \mathbb{R}$ such that $h(v) := a^{\top} v + v^{\top} B_{\mbox{\tiny N}} v$. Then, $h$ is submodular. Furthermore, inequality $\pi^{\top}v \leq T^2$ is valid for set $\{v \in \{0, 1\}^{J + J^2}: h(v) \leq T^2 \}$ for all $\pi \in \mbox{EP}_h$ and the separation of this inequality can be done by Algorithm \ref{algo:greedy}.
\end{proposition}

Note that this lifting procedure introduces $J^2$ additional variables $w_{ijk}$ for each $i$. However, $w_{ijk}$ can be treated as continuous variables when solving DCBP in view of the McCormick inequalities \eqref{mccormick}, and the number of $w_{ijk}$ variables can be reduced by half because $w_{ijk} = w_{ikj}$. In our numerical studies later, we derive more valid inequalities to strengthen the formulation in the lifted space for distributionally robust chance-constrained bin packing problems that involve DRCCs, using the bin packing structure.

{\color{black}
\begin{example}
\label{ex:polymatroidCut}
Recall Example~\ref{ex:general-cov} and the corresponding function $g(y)$ being not submodular. We set $\mu = [0,0,0]^{\top}$ and $T = 0.8$, and rewrite the constraint $g(y) \leq T$ in the form of \eqref{eq:lifted:v} as
$$0.6 v_4 + 0.2 v_6 + 0.7 v_8 + 0.1 v_9 + 0.2 v_{10} + 0.1 v_{11} + 0.6 v_{12} - 0.4 v_1 v_2 \leq 0.64,$$
where $v := [y_1, y_2, y_3, w_{11}, w_{12}, \ldots, w_{33}]^\top = [v_1,\ldots,v_{12}]^\top$. Now suppose that we are given $\hat{y} = [1, 0.5, 0.9]^{\top}$. The corresponding $\hat{v} = [1, 0.5, 0.9, 1, 0.5, 0.9, 0.5, 0.25, 0.45,$ $0.9, 0.45, 0.81]^\top$ violates the above inequality. We follow Algorithm~\ref{algo:greedy} to find an extended polymatroid inequality in the lifted space of $v$. Specifically, we sort the entries of $\hat{v}$ to obtain $\hat{v}_1 \geq \hat{v}_4 \geq \hat{v}_{3} \geq \hat{v}_6 \geq \hat{v}_{10} \geq \hat{v}_{12}\geq \hat{v}_2 \geq \hat{v}_5 \geq \hat{v}_7 \geq \hat{v}_9 \geq \hat{v}_{11} \geq \hat{v}_{8}$ and the corresponding permutation $\{1, 4, 3, 6, 10, 12, 2, 5, 7, 9, 11, 8\}$. It follows that
$\hat{\pi}_{(2)} = 0.6$, $\hat{\pi}_{(4)} = 0.2$, $\hat{\pi}_{(5)} = 0.2$, $\hat{\pi}_{(6)} = 0.6$, $\hat{\pi}_{(7)} = -0.4$, $\hat{\pi}_{(10)} = 0.1$, $\hat{\pi}_{(11)} = 0.1$, $\hat{\pi}_{(12)} = 0.7$, and all other $\hat{\pi}_{(i)}$'s equal zero. This generates the following extended polymatroid inequality
$$ 0.6 v_4 + 0.2 v_6 + 0.2 v_{10} + 0.6 v_{12}  -0.4 v_{2} + 0.1 v_{9} +0.1 v_{11} + 0.7 v_{8} \leq 0.64$$
that cuts off $\hat{v}$.

 \qedbox
\end{example}
}

\section{Numerical Studies}
\label{sec:result}
We numerically evaluate the performance of our proposed models and solution approaches. In Section \ref{sec:ACCBP}, we present the formulation of a chance-constrained bin-packing problem with DRCCs and the related 0-1 SOC reformulations. We describe the solution methods and more valid inequalities based on the bin packing structure. In Section \ref{sec:comp_setup}, we describe the experimental setup of the stochastic bin packing instances. Our results consist of two parts, which report the CPU time (Section \ref{sec:cpu}) and the out-of-sample performance of solutions given by different models (Appendix \mbox{SM3}), respectively. More specifically, Section \ref{sec:cpu} demonstrates the computational efficacy of the valid inequalities we derived in Section \ref{sec:algorithm} for the original or lifted SOC constraints. Appendix \mbox{SM3} demonstrates that DCBP solutions can well protect against the distributional ambiguity as opposed to the solutions obtained by following the Gaussian distribution assumption or by the SAA method.

{\color{black}
\subsection{Formulation of Ambiguous Chance-Constrained Bin Packing}
\label{sec:ACCBP}
For the classical bin packing problem, $[I]$ is the set of bins and $[J]$ is the set of items, where each bin $i$ has a weight capacity $T_i$ and each item $j$, if assigned to bin $i$, has a random weight $\tilde{t}_{ij}$. The deterministic bin packing problem aims to assign all $J$ items to a minimum number of bins, while respecting the capacity of each bin. If we consider a slightly more general setting by introducing a cost for each assignment, then the DCBP of bin packing with DRCCs is presented as:
\begin{subequations}\label{bin-packing}
	\begin{eqnarray}
	\label{bp:0}\min_{\mathbf z, \mathbf y} \quad & & \sum_{i=1}^I c^{\text{z}}_i z_i + \sum_{i=1}^I \sum_{j=1}^J c^{\text{y}}_{ij} y_{ij} \\
	\label{bp:1}\mbox{s.t.} \quad & & y_{ij} \leq \rho_{ij} z_i \quad \forall i \in [I], \ j \in [J] \\
	\label{bp:2} & & \sum_{i=1}^I y_{ij} = 1  \quad \forall j \in [J] \\
	\label{bp:4} & & y_{ij}, z_i \in \{0,1\}  \quad \forall i \in [I], \ j \in [J], \\
 \label{bp:3} 	& & 
	\inf_{\mathbb{P} \in \mathcal{D}} \ \mathbb{P} \Biggl\{ \sum_{j=1}^J \tilde{t}_{ij} y_{ij} \leq T_i \Biggr\} \ \geq \ 1 - \alpha_i \quad \forall i \in [I], 
	\end{eqnarray}
\end{subequations}
where $c^{\text{z}}_i$ represents the cost of opening bin $i$ and $c^{\text{y}}_{ij}$ represents the cost of assigning item $j$ to bin $i$. For each $i \in [I]$ and $j \in [J]$,  we let binary variables $z_i$ represent if bin $i$ is open (i.e., $z_i = 1$ if open and $z_i = 0$ otherwise), binary variables $y_{ij}$ represent if item $j$ is assigned to bin $i$ (i.e., $y_{ij} = 1$ if assigned and $y_{ij} = 0$ otherwise), and parameters $\rho_{ij}$ represent if we can assign item $j$ to bin $i$ (i.e., $\rho_{ij} = 1$ if we can and $\rho_{ij} = 0$ otherwise).
The objective function \eqref{bp:0} minimizes the total cost of opening bins and assigning items to bins. Constraints \eqref{bp:1} ensure that all items are assigned to open bins, constraints \eqref{bp:2} ensure that each item is assigned to one and only one bin, and constraints \eqref{bp:3} are the DRCCs.
In our computational studies, we follow Section \ref{sec:model} to derive 0-1 SOC reformulations of model \eqref{bin-packing} and then follow Section \ref{sec:algorithm} to derive valid inequalities in the original and lifted space for the 0-1 SOC reformulations.  We strengthen the extended polymatroid inequalities, as well as derive valid inequalities in the lifted space containing variables $z_i$, $y_{ij}$, and $w_{ijk}$ to further strengthen the formulation as follows. We refer to Appendices \mbox{SM1} and \mbox{SM2} for the detailed proofs of the valid inequalities below, and will test their effectiveness later.

\begin{proposition}
\label{thm:poly-z}
\begin{subequations}
For all extended polymatroid inequalities $\pi^{\top} y_i \leq T$ with regard to bin $i$, $\forall i \in [I]$, inequality
\begin{equation}
\pi^{\top} y_i \leq T z_i \label{submodular-1}
\end{equation}
is valid for the DCBP formulation. Similarly, for all extended polymatroid inequalities $\pi^{\top} v_i \leq T^2$ with regard to bin $i$, $\forall i \in [I]$, inequality
\begin{equation}
\pi^{\top} v_i \leq T^2 z_i \label{submodular-2}
\end{equation}
is valid for the DCBP formulation.
\end{subequations}
\end{proposition}

\begin{proposition}
\label{thm:poly-lifted}
Consider set
\begin{equation*}
L = \Bigl\{(z, y, w) \in \{0, 1\}^{I \times (IJ) \times (IJ^2)}: \ \mbox{\eqref{bp:1}--\eqref{bp:4}}, \ w_{ijk} = y_{ij} y_{ik}, \ \forall j, k \in [J] \Bigr\}.
\end{equation*}
Without loss of optimality, the following inequalities are valid for $L$:
\begin{subequations}
\begin{align}
w_{ijk} \ \geq & \ y_{ij} + y_{ik} + \sum_{\substack{\ell=1 \\ \ell \neq i}}^I w_{\ell jk} - 1 \quad \forall j, k \in [J] \label{lifted-1} \\
w_{ijk} \ \geq & \ y_{ij} + y_{ik} - z_i \quad \forall i \in [I], \ \forall j, k \in [J] \label{lifted-2} \\
\sum_{\substack{j=1\\ j \neq k}}^J w_{ijk} \ \leq & \ \sum_{j=1}^J y_{ij} - z_i \quad \forall i \in [I], \ \forall k \in [J] \label{lifted-3} \\
\sum_{j=1}^J \sum_{k = j+1}^J w_{ijk} \ \geq & \ \sum_{j=1}^J y_{ij} - z_i \quad \forall i \in [I]. \label{lifted-4}
\end{align}
\end{subequations}
\end{proposition}

\begin{remark}
We note that valid inequalities \eqref{lifted-1}--\eqref{lifted-4} are \emph{polynomially many} and all coefficients are in \emph{closed-form}. Hence, we do not need any separation processes for these inequalities, and we can incorporate them in the DCBP formulation without dramatically increasing its size.
\end{remark}
 }

\subsection{Computational Setup}
\label{sec:comp_setup}
We first consider $I= 6$ servers (i.e., bins) and $J = 32$ appointments (i.e., items) to test the DCBP model under various distributional assumptions and ambiguity sets. The daily operating time limit (i.e., capacity) $T_i$ of each server $i$ varies in between $[420,~540]$ minutes (i.e., 7--9 hours). We let the opening cost $c_i^{\mbox{\tiny z}}$ of each server $i$ be an increasing  function of $T_i$ such that $c_i^{\textrm{z}} = T_i^2/3600+3T_i/60$, and let all assignment costs $c_{ij}^{\textrm{y}},~\forall i\in [I],~j\in [J]$ vary in between $[0,18]$, so that the total opening cost and the total assignment cost have similar magnitudes. The above problem size and parameter settings follow the literature of surgery block allocation (see, e.g., \cite{ORA_denton2010ORA_informs, ORA_Shylo_CC, yan2013ORP}).

To generate samples of random service time (i.e., random item weight), we consider ``high mean (hM)" and ``low mean ($\ell$M)" being $25$ minutes and $12.5$ minutes, respectively. We set the coefficient of variation (i.e., the ratio of the standard deviation to the mean) as $1.0$ for the ``high variance (hV)" case and as $0.3$ for the ``low variance ($\ell$V)" case. We equally mix all four types of appointments with ``hMhV", ``hM$\ell$V", ``$\ell$MhV", ``$\ell$M$\ell$V", and thus have eight appointments of each type.
We sample 10,000 data points as the random service time of each appointment on each server, following a Gaussian distribution with the above settings of mean and standard deviation. We will hereafter call them the in-sample data. To formulate the 0-1 SOC models with diagonal covariance matrices, we use the empirical mean and standard deviation of each $\tilde{t}_{ij}$ obtained from the in-sample data and set $\alpha_i = 0.05, \ \forall i \in [I]$. Using the same $\alpha_i$-values, we formulate the 0-1 SOC models under general covariance matrices, for which we use the empirical mean and covariance matrix obtained from the in-sample data. The empirical covariance matrices we obtain have most of their off-diagonal entries being non-zero, and some being quite significant.

All the computation is performed on a Windows 7 machine with Intel(R) Core(TM) i7-2600 CPU 3.40 GHz and 8GB memory. We implement the optimization models and the branch-and-cut algorithm using commercial solver GUROBI 5.6.3 via Python 2.7.10. The GUROBI default settings are used for optimizing all the 0-1 programs, and we set the number of threads as one. When implementing the branch-and-cut algorithm, we add the violated extended polymatroid inequalities using GUROBI {\tt callback} class by \texttt{Model.cbCut()} for both integer and fractional temporary solutions. For all the nodes in the branch-and-bound tree, we generate violated cuts at each node as long as any exists. The optimality gap tolerance is set as $0.01\%$. We also set the threshold for identifying violated cuts as $10^{-4}$, and set the time limit for computing each instance as 3600 seconds.

\subsection{CPU Time Comparison}
\label{sec:cpu}

We solve 0-1 SOC reformulations or approximations of DCBP, and use either a diagonal or a general covariance matrix in each model. Our valid inequalities significantly reduce the solution time of directly solving the 0-1 SOC models in GUROBI, while the extended polymatroid inequalities generated based on the approximate and lifted SOC constraints perform differently depending on the problem size. The details are presented as follows.

\subsubsection{Computing 0-1 SOC models with diagonal matrices}
We first optimize 0-1 SOC models with a diagonal matrix in constraint \eqref{soc}, of which the left-hand side function $g(y)$ is submodular, and thus we use extended polymatroid inequalities \eqref{submodular-1} with $\pi \in EP_g$ in a branch-and-cut algorithm. Table~\ref{tab1:diag_cpu} presents the CPU time (in seconds), optimal objective values, and other solution details (including ``{\bf Server}" as the number of open servers, ``{\bf Node}" as the total number of branching nodes,  and ``{\bf Cut}" as the total number of cuts added) for solving the three 0-1 SOC models {\bf DCBP1} (i.e., 
using ambiguity set $\mathcal D_1$), {\bf DCBP2} (i.e., 
using ambiguity set $\mathcal D_2$ with $\gamma_1 = 1, \ \gamma_2 = 2$), and {\bf Gaussian} (assuming Gaussian distributed service time). We also implement the SAA approach (i.e., row ``{\bf SAA}'') by optimizing the MILP reformulation of the chance-constrained bin packing model based on the 10,000 in-sample data points.
We compare the branch-and-cut algorithm using our extended polymatroid inequalities (in rows ``{\bf B\&C}'') with directly solving the 0-1 SOC models in GUROBI (in rows ``{\bf w/o Cuts}'').

\begin{table}[htbp]
	\centering
	\caption{CPU time and solution details for solving instances with diagonal matrices}
    \resizebox{\textwidth}{!}{%
	\begin{tabular}{ccrrrrrr}
		\hline
		Approach & Model & \multicolumn{1}{c}{Time (s)} & \multicolumn{1}{c}{Opt.\ Obj.\ } & \multicolumn{1}{c}{Server}& \multicolumn{1}{c}{Opt.\ Gap} & \multicolumn{1}{c}{Node} & \multicolumn{1}{c}{Cut}\\
		\hline
		\multirow{3}[1]{*}{B\&C} & \multicolumn{1}{c}{DCBP1} & 0.73  & 328.99 &3& 0.00\% & 83    & 82 \\
		& \multicolumn{1}{c}{DCBP2} & 27.50 & 366.54&3 & 0.00\% & 2146  & 2624 \\
		& \multicolumn{1}{c}{Gaussian} & 0.13  & 297.94&2 & 0.00\% & 0     & 0 \\
		\hline
		\multirow{3}[1]{*}{w/o Cuts} & \multicolumn{1}{c}{DCBP1} & 95.73 & 328.99 &3& 0.01\% & 76237 & N/A \\
		& \multicolumn{1}{c}{DCBP2} & LIMIT & 380.09&2 & 9.15\% & 409422 & N/A \\
		& \multicolumn{1}{c}{Gaussian} & 0.02  & 297.94&2 & 0.00\% & 16    & N/A \\
		\hline
		SAA  & \multicolumn{1}{c}{MILP} & 21.20 & 297.94 &2& 0.00\% & 89    & N/A \\
		\hline
	\end{tabular}%
    }
	\label{tab1:diag_cpu}%
\end{table}%

In Table~\ref{tab1:diag_cpu}, the branch-and-cut algorithm quickly optimizes DCBP1 and DCBP2. Especially,  if being directly  solved by GUROBI, DCBP2 cannot be solved within the 3600-second time limit and ends with a $9.15\%$ optimality gap. Solving DCBP1 by using the branch-and-cut algorithm is much faster than solving the large-scale SAA-based MILP model, while the solution time of DCBP2 is similar to the latter. The two DCBP models also yield higher objective values, since they both provide more conservative solutions that open one more server than either the Gaussian or the SAA-based approach.

\subsubsection{Computing 0-1 SOC models with general covariance matrices}
In this section, we focus on testing DCBP2 yielded by the ambiguity set $\mathcal D_2$ with parameters  $\gamma_1=1$, $\gamma_2=2$, and $\alpha_i = 0.05, \ \forall i\in [I]$. We use empirical covariance matrices of the in-sample data. Note that these covariance matrices are general and non-diagonal. We compare the time of solving the 0-1 SOC reformulations of DCBP2 on ten independently generated instances. We examine two implementations of the branch-and-cut algorithm: one uses extended polymatroid inequalities \eqref{submodular-1} with $\pi \in EP_{g^L}$ based on the relaxed 0-1 SOC constraint \eqref{soc-relaxed}, and the other uses the extended polymatroid inequalities \eqref{submodular-2} based on the lifted SOC constraint.

Table~\ref{tab2:general_cpu} reports the CPU time (in seconds) and number of branching nodes (in column ``\textbf{Node}") for various methods. First, we directly solve the 0-1 SOC models of DCPB2 in GUROBI, without or with the linear valid inequalities \eqref{lifted-1}--\eqref{lifted-4}, and report their results in columns ``\textbf{w/o Cuts}'' and ``{\bf Ineq.}'', respectively. We then implement the branch-and-cut algorithm, and examine the results of using extended polymatroid inequalities \eqref{submodular-1} with $\pi \in EP_{g^L}$ (reported in columns ``{\bf B\&C-Relax}'') and cuts \eqref{submodular-2} based on lifted SOC constraints (reported in columns ``{\bf B\&C-Lifted}''). {\color{black} (The latter does not involve valid inequalities \eqref{submodular-1} used in the former.) The time reported under {\bf B\&C-Relax} also involves the time of solving SDPs for obtaining  $\Delta^{\mbox{\tiny L}}$ and the relaxed 0-1 SOC constraint \eqref{soc-relaxed}. The time of solving related SDPs are small (varying from 1 to 2 seconds for instances of different sizes) and negligible as compared to the total B\&C time.}
For both B\&C methods, we also present the number of extended polymatroid inequalities (see column ``\textbf{Cut}") added.

\begin{table}[htbp]
	\centering
	\caption{CPU time of DCBP2 solved by different methods with general covariance matrices}
    \resizebox{\textwidth}{!}{%
	\begin{tabular}{crr|rr|rrr|rrr}
		    	\hline
		    	\multicolumn{1}{c}{\multirow{2}[2]{*}{Instance}} & \multicolumn{2}{c|}{w/o Cuts} & \multicolumn{2}{c|}{Ineq.} &\multicolumn{3}{c|}{B\&C-Relax}  & \multicolumn{3}{c}{B\&C-Lifted}   \\
		    	\multicolumn{1}{c}{} & Time (s) & Node  & Time (s) &  Node  & Time (s) &  Node  &  Cut & Time (s) & Node  & Cut \\
		    	\hline
		    	1     & 286.29 & 10409 & 156.50 & 795    & 51.99 & 9095 & 702 & {\bf 35.03} & 618   & 823   \\
		    	2     & 433.32 & 10336 & 167.91 & 687  & 26.63 & 6524   & 698  & {\bf 12.34} & 405   & 235    \\
		    	3     & 284.17 & 10434 & 206.82 & 971    & 70.43 & 17420   & 621 & {\bf 29.84} & 595   & 729   \\
		    	4     & 310.11 & 10302 & 139.06 & 656   & {\bf 15.37} & 2467   & 723 & 25.31 & 419   & 617    \\
		    	5     & 329.32 & 10453 & 181.83 & 777   & 56.53 & 12349  & 737 & {\bf 35.09} & 678   & 921    \\
		    	6     & 365.28 & 10300 & 168.26 & 652   & {\bf 23.89} & 4807   & 695 & 26.73 & 555   & 595    \\
		    	7     & 296.55 & 10759 & 198.87 & 873    & 45.21 & 11585   & 738  & {\bf 21.08} & 440   & 626 \\
		    	8     & 278.62 & 10490 & 211.05 & 900   & 53.84 & 14540  & 721 & {\bf 47.78} & 1064  & 1686  \\
		    	9     & 139.24 & 7771  & 177.41 & 632     & 19.90 & 3918  & 645 & {\bf 19.37} & 216   & 360  \\
		    	10    & 297.72 & 10330 & 159.52 & 822  & 30.36 & 6877  & 649 & {\bf 29.43} & 400   & 727    \\
		    	\hline
	\end{tabular}%
    }
	\label{tab2:general_cpu}%
\end{table}%

In Table~\ref{tab2:general_cpu}, we highlight the solution time of the method that runs the fastest for each instance. Note that without the extended polymatroid inequalities or the valid inequalities, the GUROBI solver takes the longest time for solving all the instances except instance \#9. Adding the valid inequalities \eqref{lifted-1}--\eqref{lifted-4} to the solver reduces the solution time by 40\% or more in almost all the instances, and drastically reduces the number of branching nodes. The extended polymatroid inequalities \eqref{submodular-1}--\eqref{submodular-2} further reduce the CPU time significantly (see columns ``{\bf B\&C-Relax}'' and ``{\bf B\&C-Lifted}''). Moreover, for all the instances having 6 servers and 32 appointments, the algorithm using the extended polymatroid inequalities \eqref{submodular-2} runs faster in eight out of ten instances than the algorithm using cuts \eqref{submodular-1} based on relaxed SOC constraints without lifting. It indicates that the extended polymatroid inequalities generated by the lifted SOC constraints are more effective than those generated by the relaxed SOC constraints. This observation is overturned when we later increase the problem size.

In the following, we continue reporting the CPU time of solving DCBP2 with general covariance matrix. We vary the problem sizes (i.e., values of $I$ and $J$) in Section \ref{sec:time:size}, and vary the values of $\Lambda$ in the SOC constraint \eqref{soc} in Section \ref{sec:time:lambda}.

\subsubsection{Solving 0-1 SOC models under different problem sizes}
\label{sec:time:size}

We use the same problem settings as in Section \ref{sec:comp_setup}, and vary $I = 6, 8, 10$ and $J= 32, 40$ to test DCBP2 instances with different sizes. We still keep an equal mixture of all the four appointment types in each instance. Table~\ref{tab3:size} presents the computational time (in seconds), the total number of branching nodes (``\textbf{Node}"), and the total number of extended polymatroid inequalities generated (``\textbf{Cuts}"; if applicable) for solving the 0-1 SOC reformulation of DCBP2 by directly using GUROBI (``\textbf{w/o Cuts}'') and by using the two implementations of the extended polymatroid inequalities (``\textbf{B\&C-Relax}'' and ``\textbf{B\&C-Lifted}'').

\begin{table}[htbp]
  \centering
  \caption{CPU time of DCBP2 with general covariance matrices for different problem sizes}
  \resizebox{\textwidth}{!}{%
    \begin{tabular}{cccrrrrrrrrrr}
    \hline
          &    \multirow{2}{*}{Method}   &       & \multicolumn{5}{c}{$J=32$}  & \multicolumn{5}{c}{$J=40$}  \\
      \cline{4-13}
          &       & Inst.\ & 1     & 2     & 3     & 4     & 5     & 6     & 7     & 8     & 9     & 10 \\
          \hline
    \multirow{8}{*}{$I=6$}  & \multirow{3}{*}{B\&C-Relax}  & Time (s) & 51.99	 & 26.63	&70.43	& {\bf 15.37}	&56.53	&{\bf 6.87}	&{\bf 12.76} &	{\bf 1.59}	&{\bf 2.36}&	{\bf 12.73}  \\
          &       & Node  & 9095 &	6524&	17420	&2467	&12349&	1009	&1322	&176&	285&	1270 \\
          &       & Cut   & 702&	698	&621 &	723	&737&	174&	604&	171	&179&	602 \\
             \cline{3-13}

         &     \multirow{3}{*}{B\&C-Lifted}  & Time (s) & {\bf 35.03} & {\bf 12.34} & {\bf 29.84} & 25.31 & {\bf 35.09} & 64.58 & 98.18 & 91.12 & 60.11 & 59.50 \\
          &       & Node  & 618 &405& 595 &419 &678 &274 &484 &447& 289& 234 \\
          &       & Cut   & 823&  235 &  729 &  617  & 921 & 470 & 690 & 688 & 462 & 394\\
                       \cline{3-13}

          & \multirow{2}{*}{w/o Cuts} & Time (s) & 286.29 & 433.32 & 284.17 & 310.11 & 329.32 & 1654.31 & 208.12 & 1182.46 & 1580.41 & 1266.27 \\
          &       & Node  & 10409 & 10336 & 10434 & 10302 & 10453 & 10525 & 1272  & 10732 & 10658 & 10642 \\
          \hline
    \multirow{8}{*}{$I=8$}   & \multirow{3}{*}{B\&C-Relax}  & Time (s) & {\bf 41.57}	&139.41&	{\bf 55.22} &	261.24&	305.72	&{\bf 23.91}	&{\bf 9.73}	&{\bf 17.76}	&{\bf 27.16}&	 {\bf 12.98}\\
          &       & Node  & 8342	&29042	&12267	&49820	&61334	&2130&	1240&	1561&	2607&	1024 \\
          &       & Cut   & 737	&770&	742	&803&	790	&714&	199&	728	&702	&690\\
                       \cline{3-13}
          &              \multirow{3}{*}{B\&C-Lifted}  & Time (s) & 106.03 &{\bf 28.55} & 84.64 &{\bf 97.05} & {\bf 13.56} & 331.29 & 273.14 &307.06 &178.41 &161.39 \\
          &       & Node  & 678 & 502 & 647 & 634  &125  &1177  &836  &1397  &457  &529 \\
          &       & Cut   & 114 & 691&  128 & 143 & 216 & 1781 & 1175&  2066 & 703 & 719\\
                       \cline{3-13}

          & \multirow{2}{*}{w/o Cuts}  & Time (s) & 866.12 & 597.43 & 649.72 & 683.18 & 497.15 & 2265.53 & 2428.60 & 2294.62 & 1781.95 & 851.99 \\
          &       & Node  & 10338 & 10305 & 10309 & 10306 & 14386 & 11441 & 11219 & 11708 & 11128 & 5241 \\
                    \hline
    \multirow{8}{*}{$I=10$} & \multirow{3}{*}{B\&C-Relax}  & Time (s) & {\bf 3.75}	& {\bf 9.28}& {\bf 	6.56}	& {\bf 3.23}& 	{\bf 16.71} & 	{\bf 29.94}	& {\bf 80.34}	& {\bf 22.58}	 & {\bf 24.48}	& {\bf 339.93} \\
          &       & Node  & 637	&972&	659	&549&	2274&	2336	&7315&	1870&	1959	& 34306 \\
          &       & Cut   & 241	 &552 &	390	 &230  &	741	 &767	 &714 &	736 &	729	 &715\\
                       \cline{3-13}

            &            \multirow{3}{*}{B\&C-Lifted}  & Time (s) & 108.43 & 117.44 & 120.60 & 22.10 & 111.37 & 186.72 & 714.45 & 197.42 & 549.90 & 661.13 \\
          &       & Node  & 668   & 785   & 828   & 291   & 779   & 766   & 1108  & 811   & 896   & 1209 \\
          &       & Cut   & 108   & 191   & 314   & 281   & 188   & 1196  & 850   & 1106  & 568   & 808 \\
                       \cline{3-13}

          & \multirow{2}{*}{w/o Cuts}  & Time (s) & 987.92 & 1140.23 & 183.06 & 1113.09 & 1425.83 & 2382.97 & 2917.03 & LIMIT & 2052.42 & 2451.62 \\
          &       & Node  & 10353 & 10357 & 4992  & 10307 & 10401 & 11015 & 11197 & 12101 & 10812 & 11001 \\
                    \hline
    \end{tabular}%
    }
  \label{tab3:size}%
\end{table}%

In Table~\ref{tab3:size}, we again highlight the solution time of the method that runs the fastest in each instance. We keep the first five instances we reported in Table~\ref{tab2:general_cpu} for instances with $I=6, \ J =32$, and report five instances for other $(I, J)$ combinations. From Table~\ref{tab3:size}, we observe that both implementations of the extended polymatroid inequalities run significantly faster than directly using GUROBI, especially when we increase the problem sizes (i.e., $I$ increased from 6 to 10, and $J$ increased from 32 to 40). In particular, the CPU time of directly using GUROBI is consistently 1 or 2 orders of magnitude larger than that of our approaches. For smaller $(I, J)$-values (e.g., $(I, J) = (6, 32)$ or $(I, J) =(8, 32)$), we see that {\bf B\&C-Lifted} sometimes runs faster than {\bf B\&C-Relax}, but for all the other $(I, J)$ combinations, the latter completely dominates the former. This is expected because the cuts \eqref{submodular-2} are generated in a lifted space with $J^2$ additional variables for each server $i \in [I]$. Therefore, it makes sense that the scalability of {\bf B\&C-Lifted} is worse than that of {\bf B\&C-Relax}, which uses cuts \eqref{submodular-1} without lifting.

\subsubsection{Solving 0-1 SOC models with different $\Lambda$-values in \eqref{soc}}
\label{sec:time:lambda}
We again focus on instances with $I = 6$ and $J = 32$ under the same general covariance matrix $\Sigma$ obtained from the in-sample data points. We let $\Lambda := \Omega^2 \Sigma$ in the 0-1 SOC constraint \eqref{soc} and adjust the scalar $\Omega$ to obtain different $\Lambda$. We want to show how the computational time of directly using GUROBI increases as we increase $\Omega$, as compared to using the branch-and-cut algorithm with extended polymatroid inequalities. (The results of {\bf B\&C-Lifted} are used here and similar observations can be made if the results of {\bf B\&C-Relax} are used.)

Considering specific cases of the SOC constraint \eqref{soc} for modeling DCBP, we have $\Omega= \Phi^{-1}(1-\alpha)= 1.64$ for the Gaussian approximation model when $\alpha = 0.05$, and $\Omega=\sqrt{\gamma_2/\alpha} = 6.32$ for DCBP2 when $\gamma_2 = 2$ and $\alpha = 0.05$. We test four values of $\Omega$ equally distributed in between $[1.64,6.32]$ including the two end points. 
	\begin{figure}[htbp]
		\centering
      \includegraphics[width =0.8\textwidth]{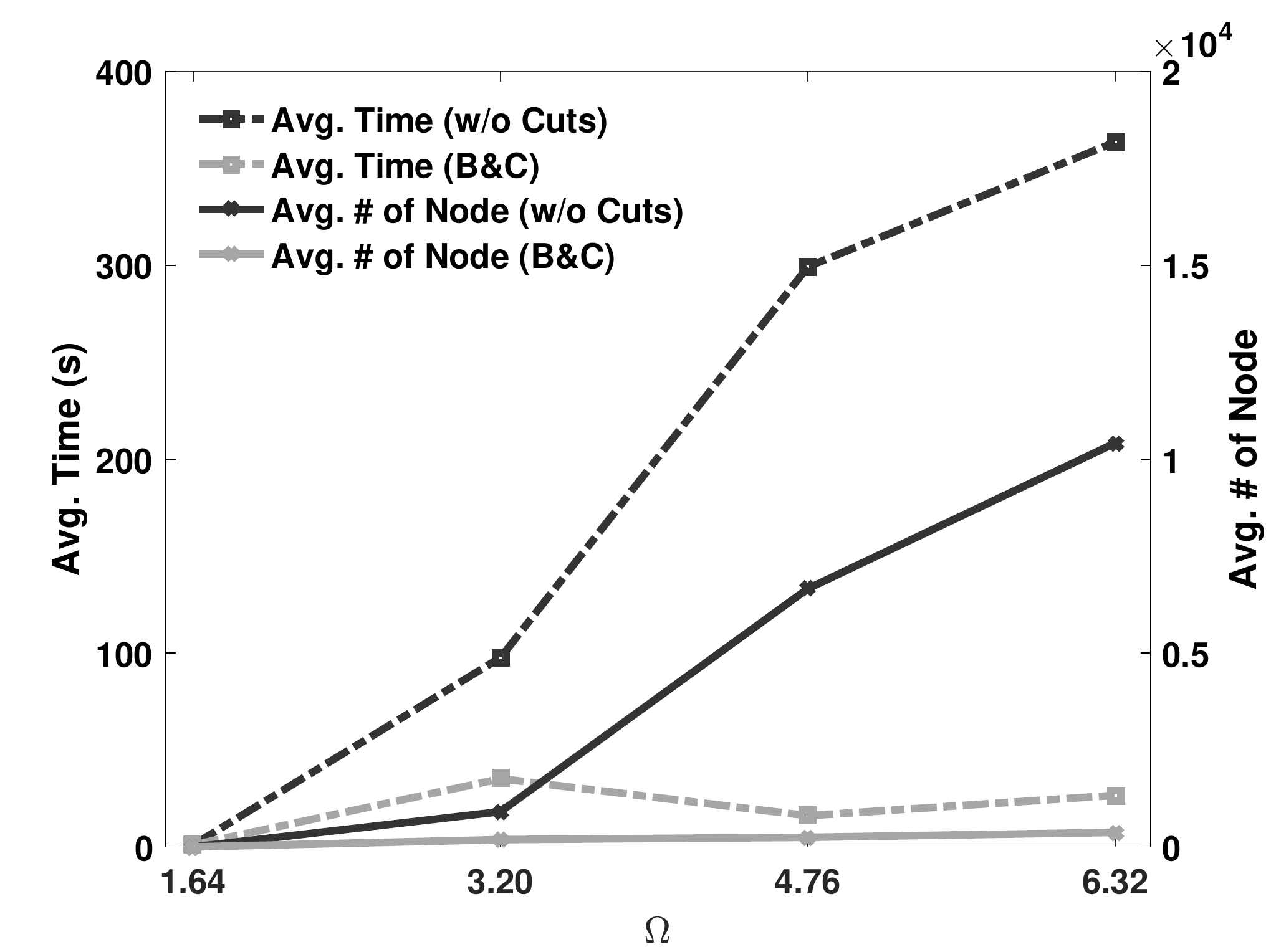}
		\caption{CPU time and number of branching nodes under different $\Omega$-values in constraint \eqref{soc}}
		\label{fig1:omega}
	\end{figure}
	
Figure~\ref{fig1:omega} depicts the average CPU time and the average number of branching nodes of solving five independent DCBP2 instances for each $\Omega$-setting.
Specifically, the four values $1.64$, $3.20$, $4.76$, and $6.32$ of $\Omega$ lead to $0.98$, $97.68$, $299.29$, and $363.56$ CPU seconds when directly using GUROBI, respectively, together with the significantly growing number of branching nodes $0$, $909.8$, $6662$, and $10418$, respectively. On the other hand, the branch-and-cut algorithm with the extended polymatroid inequalities respectively takes $1.03$, $35.19$, $16.09$, and $26.63$ CPU seconds on average for solving the same instances, and branches on average $0$, $188.6$, $247.2$, and $374.6$ nodes, respectively. This indicates that our approach is more scalable than directly using the off-the-shelf solvers.

{\color{black}
\section{Conclusions}
\label{conclusion}

In this paper, we considered distributionally robust individual chance constraints, where the true distributional information of the constraint coefficients is ambiguous and only the empirical first and second moments are given. The goal is to restrict the worst-case probability of violating a linear constraint under a given threshold. We provided 0-1 SOC representations of DRCCs under two types of ambiguity sets. In addition, we derived an efficient way of obtaining extended polymatroid inequalities for the 0-1 SOC constraints in both original and lifted spaces. Via extensive numerical studies, we demonstrated that our solution approaches significantly accelerate solving the DCBP model as compared to the state-of-the-art commercial solvers.  
In particular, a branch-and-cut algorithm with extended polymatroid inequalities in the original space scales very well as the problem size grows.

For future research, 
we plan to investigate DRCCs under other types of ambiguity sets, which could take into account not only the moment information but also density or structural information of their probability distribution. The connections between SOC program, SDP, and submodular optimization are also interesting to study.
}


\section*{Acknowledgments}
This work is supported in part by the National Science Foundation under grants CMMI-1433066 and CMMI-1662774. The authors are grateful to Profs.\ Alper Atamt{\"u}rk and Andr{\'e}s G{\'o}mez for their comments on an earlier version of this paper and for bringing to our attention the works~\cite{atamturk2018network,atamturk2017polymatroid,bhardwaj2015binary}.
The authors are grateful to the two referees and the Associate Editor for their constructive comments and helpful suggestions.

\newpage

\section*{Supplementary Materials}

\noindent {\bf SM1. Proof of Proposition \ref{thm:poly-z}} \label{sec:proofs} \\

\beginproof When $z_i = 1$, inequality \eqref{submodular-1} reduces to the extended polymatroid inequality. When $z_i = 0$, we have $y_{ij} = 0$ for all $j \in [J]$ due to constraints \eqref{bp:1}. It follows that inequality \eqref{submodular-1} holds valid.

When $z_i = 1$, inequality \eqref{submodular-2} reduces to the extended polymatroid inequality. When $z_i = 0$, we have $y_{ij} = 0$ for all $j \in [J]$ due to constraints \eqref{bp:1} and so $w_{ijk} = 0$ for all $j, k \in [J]$. It follows that $v_i = 0$ by definition. Hence, inequality \eqref{submodular-2} holds valid. \qedbox

\noindent {\bf SM2. Proof of Proposition \ref{thm:poly-lifted}} \label{sec:proofs2} \\

\beginproof ({\bf Validity of inequality \eqref{lifted-1}}) If $j = k$, then $w_{ijk} = y_{ij}^2 = y_{ij}$. In this case, inequality \eqref{lifted-1} reduces to $y_{ij} \geq 2y_{ij} + \sum_{\substack{\ell=1 \\ \ell \neq i}}^I y_{\ell j} - 1$, which clearly holds because $y_{ij} + \sum_{\substack{\ell=1 \\ \ell \neq i}}^I y_{\ell j} = \sum_{i=1}^I y_{ij} \leq 1$. If $j \neq k$, then we discuss the following two cases:
\begin{enumerate}[1.]
\item If $\max\{y_{ij}, y_{ik}\} = 1$, then we assume $y_{ij} = 1$ without loss of generality. It follows that $y_{\ell j} = 0$ due to constraints \eqref{bp:2} and so $w_{\ell jk} = y_{\ell j} y_{\ell k} = 0$ for all $\ell \neq i$. Hence, $\sum_{\substack{\ell=1 \\ \ell \neq i}}^I w_{\ell jk} = 0$ and inequality \eqref{lifted-1} reduces to $w_{ijk} \geq y_{ij} + y_{ik} - 1$, which holds valid.
\item If $\max\{y_{ij}, y_{ik}\} = 0$, then $y_{ij} = y_{ik} = 0$ and $w_{ijk} = 0$. It remains to show $\sum_{\substack{\ell=1 \\ \ell \neq i}}^I w_{\ell jk} \leq 1$. Indeed, since $w_{\ell jk} \leq y_{\ell j}$, we have $\sum_{\substack{\ell=1 \\ \ell \neq i}}^I w_{\ell jk} \leq \sum_{\substack{\ell=1 \\ \ell \neq i}}^I y_{\ell j} \leq \sum_{\ell=1}^I y_{\ell j} = 1$, where the last equality is due to constraints \eqref{bp:2}.
\end{enumerate}

\noindent({\bf Validity of inequality \eqref{lifted-2}}) This inequality clearly holds valid when $z_i = 1$. When $z_i = 0$, we have $y_{ij} = y_{ik} = 0$ due to constraints \eqref{bp:1}. It follows that $w_{ijk} = y_{ij} y_{ik} = 0$ and so the inequality holds valid.

\noindent({\bf Validity of inequality \eqref{lifted-3}}) This inequality holds valid when $z_i = 0$. When $z_i = 1$, this inequality is equivalent to $y_{ik}\sum_{\substack{j=1\\ j \neq k}}^J y_{ij} \leq \sum_{j=1}^J y_{ij} - 1$. We discuss the following two cases:
\begin{enumerate}[1.]
\item If $y_{ik} = 0$, then $\sum_{j=1}^J y_{ij} \geq 1$ without loss of optimality because $z_i = 1$. Inequality \eqref{lifted-3} holds valid.
\item If $y_{ik} = 1$, then $y_{ik}\sum_{\substack{j=1\\ j \neq k}}^J y_{ij} = \sum_{\substack{j=1\\ j \neq k}}^J y_{ij}$. Meanwhile, $\sum_{j=1}^J y_{ij} - 1 = \sum_{\substack{j=1\\ j \neq k}}^J y_{ij} + y_{ik} - 1 = \sum_{\substack{j=1\\ j \neq k}}^J y_{ij}$. Inequality \eqref{lifted-3} holds valid.
\end{enumerate}

\noindent({\bf Validity of inequality \eqref{lifted-4}}) This inequality holds valid when $z_i = 0$. When $z_i = 1$, we have $\sum_{j=1}^J y_{ij} \geq 1$ without loss of optimality. It follows that $\sum_{j=1}^J y_{ij} = 1$ or $\sum_{j=1}^J y_{ij} \geq 2$, and so $(\sum_{j=1}^J y_{ij} - 1)(\sum_{j=1}^J y_{ij} - 2) \geq 0$. Hence,
\begin{align*}
\Bigl(\sum_{j=1}^J y_{ij} - 1\Bigr)\Bigl(\sum_{j=1}^J y_{ij} - 2\Bigr) \ = & \ \sum_{j=1}^J y^2_{ij} + 2\sum_{j=1}^J \sum_{{\color{black}k = j+1}}^J y_{ij} y_{ik} - 3\Bigl(\sum_{j=1}^J y_{ij}\Bigr) + 2 \\
= & \ \sum_{j=1}^J y_{ij} + 2\sum_{j=1}^J \sum_{{\color{black}k = j+1}}^J w_{ijk} - 3\Bigl(\sum_{j=1}^J y_{ij}\Bigr) + 2 \\
= & \ 2\sum_{j=1}^J \sum_{{\color{black}k = j+1}}^J w_{ijk} - 2\biggl[\Bigl(\sum_{j=1}^J y_{ij}\Bigr) - 1\biggr] \ \geq \ 0.
\end{align*}
Inequality \eqref{lifted-4} follows. \qedbox

\noindent {\bf SM3. Out-of-Sample Performance of DCBP} \label{sec:out-of-sample}\\

Through testing instances of chance-constrained bin packing, we show that DCBP solutions have very low probabilities of violating capacities in all the out-of-sample tests, even when the distributional information is misspecified. Specifically, we evaluate the out-of-sample performance of the optimal solutions to the DCBP1, DCBP2, Gaussian-based 0-1 SOC models, and the SAA-based MILP model.
To generate the out-of-sample reference scenarios, we consider either misspecified distribution type or misspecified moment information as follows.
\begin{itemize}
\item {\bf Misspecified distribution type:} We sample 10,000 out-of-sample data points from a two-point distribution having the same mean and standard deviation of each random variable $\tilde{t}_{ij}$ for $i\in [I]$ and $j\in [J]$ as the in-sample data. The service time is realized as $\mu_{ij} + \frac{(1-p)}{\sqrt{p(1-p)}} \sigma_{ij}$ with probability $p$ ($0 < p < 1$) and as $\mu_{ij} - \frac{\sqrt{p(1-p)}}{(1-p)} \sigma_{ij}$ with probability $1-p$, where $\mu_{ij}$ and $\sigma_{ij}$ are the sample mean and standard deviation of $\tilde{t}_{ij}$ obtained from the in-sample data. We set $p = 0.3$ so that we have smaller probability of having larger service time realizations.
\item {\bf Misspecified moments:} Alternatively, we sample 10,000 data points from the Gaussian distribution, but only consider the hM$\ell$V type of appointments, instead of an equal mixture of all the four types. In each sample and for each $i \in [I]$, we draw a standard-Gaussian random number $\rho_i$, and for each $j\in [J]$, generate a service time realization as $\mu_{ij} + \rho_i\sigma_{ij}$.
\end{itemize}

\paragraph{Performance of solutions under diagonal matrices}
Under diagonal matrices, both of the two DCBP models open three servers (i.e., Servers 4, 5, 6 by DCBP1 and Servers 2, 4, 5 by DCBP2), while the Gaussian and SAA approaches only open Servers 4 and~6. We first use the $10,000$ out-of-sample data points given by misspecified distribution type, namely, the two-point distribution. Table~\ref{tab4:discrete} reports each solution's probability of having the total time of assigned appointments not exceeding the capacity of the server to which they are assigned.
%
\begin{table}[htbp]
	\centering
	\caption{Solution reliability in out-of-sample data following a misspecified distribution type}
	\begin{tabular}{crrrr}
		\hline
		Model & \multicolumn{1}{c}{Server 2} & \multicolumn{1}{c}{Server 4} & \multicolumn{1}{c}{Server 5} & \multicolumn{1}{c}{Server 6} \\
		\hline
		DCBP1 & N/A    & 1.00 & 1.00  & 1.00 \\
		DCBP2 & 1.00  & 1.00  & 1.00  & N/A \\
		Gaussian & N/A    & 0.69  & N/A    & 0.91 \\
		SAA   & N/A    & 0.69  & N/A    & 0.91 \\
		\hline
	\end{tabular}%
	\label{tab4:discrete}%
	~\\``N/A'': the server is not opened by using the corresponding method.
\end{table}%

Recall that $\alpha_i=0.05$ for all $i$ used in all four approaches. The reliability results of the Gaussian and SAA approaches are significantly lower than the desired probability threshold $1-\alpha_i=0.95$ on Server 4, and slightly lower than $0.95$ on Server 6. On the other hand, the optimal solutions of DCBP1 and DCBP2 do not exceed the capacity of any open servers.

Next, we use the $10,000$ out-of-sample data points given by misspecified moments. Table~\ref{tab5:moment_amb} reports the reliability performance of each optimal solution.
\begin{table}[htbp]
	\centering
	\caption{Solution reliability in out-of-sample scenarios with misspecified moments}
	\begin{tabular}{crrrr}
		\hline
		Model & \multicolumn{1}{c}{Server 2} & \multicolumn{1}{c}{Server 4} & \multicolumn{1}{c}{Server 5} & \multicolumn{1}{c}{Server 6} \\
		\hline
		DCBP1 & N/A    & 0.94  & 1.00  & 1.00 \\
		DCBP2 & 0.98  & 1.00  & 0.99  & N/A \\
		Gaussian & N/A    & 0.59  & N/A    & 0.89 \\
		SAA   & N/A    & 0.59  & N/A    & 0.89 \\
		\hline
	\end{tabular}%
	\label{tab5:moment_amb}%
	~\\``N/A'': the server is not opened by using the corresponding method.
\end{table}%
The DCBP2 solution still outperforms solutions given by all the other approaches and achieves the desired reliability in all the three open servers. The Gaussian and SAA solutions perform poorly when the moment information is different from the empirical inputs. The DCBP1 solution respects the capacities of Servers 5 and 6 with sufficiently high probability (i.e., $> 0.95$), but yields a slightly lower reliability ($0.94$) than the threshold on Server 4.

\paragraph{Performance of solutions under general matrices}
We optimize all the models under general matrices by using the empirical covariance matrices of the in-sample data, and report their corresponding solutions in Table~\ref{tab:solution:general}. Each entry illustrates the number of appointments assigned to an open server. Note that the Gaussian and SAA approaches yield the same solution of opening servers and assigning appointments.
\begin{table}[htbp]
	\centering
	\caption{Optimal open servers and appointment-to-server assignments under general matrices}
	\begin{tabular}{ccccc}
		\hline
		Model & \multicolumn{1}{c}{Server 3} & \multicolumn{1}{c}{Server 4} & \multicolumn{1}{c}{Server 5} & \multicolumn{1}{c}{Server 6} \\
		\hline
		DCBP1 & 12   & N/A  & 13  & 7 \\
		DCBP2 & 12   & 11 & 9  & N/A \\
		Gaussian & 15 & N/A    & 17  & N/A \\
		SAA   &15 & N/A    & 17  & N/A \\
		\hline
	\end{tabular}%
	\label{tab:solution:general}%
    ~\\``N/A'': the server is not opened by using the corresponding method.
\end{table}%

We test the solutions shown in Table~\ref{tab:solution:general} in the out-of-sample scenarios under misspecified distribution type, and present their reliability performance in Table~\ref{tab:out-of-sample:general}. We again show that under general covariance matrices, the DCBP2 model yields the most conservative solution that does not exceed any open server's capacity, while DCBP1 only ensures the desired reliability on Servers 3 and 6, but not on Server 5. The Gaussian and SAA approaches cannot produce solutions that can achieve the desired reliability threshold on any of their open servers.

\begin{table}[htbp]
	\centering
	\caption{Solution reliability in out-of-sample scenarios with misspecified moments}
	\begin{tabular}{crrrr}
		\hline
		Model & \multicolumn{1}{c}{Server 3} & \multicolumn{1}{c}{Server 4} & \multicolumn{1}{c}{Server 5} & \multicolumn{1}{c}{Server 6} \\
		\hline
		DCBP1 & 1.00  & N/A  & 0.91  & 1.00 \\
		DCBP2 & 1.00 & 1.00  & 1.00  & N/A \\
		Gaussian & 0.91 & N/A  & 0.91    & N/A \\
		SAA   & 0.91 & N/A  & 0.91    & N/A \\
		\hline
	\end{tabular}%
	\label{tab:out-of-sample:general}%
	~\\``N/A'': the server is not opened by using the corresponding method.
\end{table}%

\newpage

\bibliographystyle{siamplain}

\end{document}